

\overfullrule=0pt


\font\headerfont=cmr10

\font\smallfont=cmr7
\font\sectionfont=cmbx12

\newlinechar=`@
\def\forwardmsg#1#2#3{\immediate\write16{@*!*!*!* forward reference should
be: @\noexpand\forward{#1}{#2}{#3}@}}
\def\nodefmsg#1{\immediate\write16{@*!*!*!* #1 is an undefined reference@}}

\def\forwardsub#1#2{\def\newref{{#2}{#1}}}

\def\forward#1#2#3{%
\expandafter\expandafter\expandafter\forwardsub\expandafter{#3}{#2}
\expandafter\ifx\csname#1\endcsname\relax\else%
\expandafter\ifx\csname#1\endcsname\newref\else%
\expandafter\ifx\csname#2\endcsname\relax\else%
\forwardmsg{#1}{#2}{#3}\fi\fi\fi%
\expandafter\let\csname#1\endcsname\newref}

\def\firstarg#1{\expandafter\argone #1}\def\argone#1#2{#1}
\def\secondarg#1{\expandafter\argtwo #1}\def\argtwo#1#2{#2}

\def\ref#1{\expandafter\ifx\csname#1\endcsname\relax%
  {\nodefmsg{#1}\bf`#1'}\else%
  \expandafter\firstarg\csname#1\endcsname%
  ~\htmllocref{#1}{\expandafter\secondarg\csname#1\endcsname}\fi}

\def\refscor#1{\expandafter\ifx\csname#1\endcsname\relax
  {\nodefmsg{#1}\bf`#1'}\else
  Corollaries~\htmllocref{#1}{\expandafter\secondarg\csname#1\endcsname}\fi}

\def\refs#1{\expandafter\ifx\csname#1\endcsname\relax
  {\nodefmsg{#1}\bf`#1'}\else
  \expandafter\firstarg\csname #1\endcsname
  s~\htmllocref{#1}{\expandafter\secondarg\csname#1\endcsname}\fi}

\def\refn#1{\expandafter\ifx\csname#1\endcsname\relax
  {\nodefmsg{#1}\bf`#1'}\else
  \htmllocref{#1}{\expandafter\secondarg\csname #1\endcsname}\fi}

\def\pageref#1{\expandafter\ifx\csname#1\endcsname\relax
  {\nodefmsg{#1}\bf`#1'}\else
  \expandafter\firstarg\csname#1\endcsname%
  ~\htmllocref{#1}{\expandafter\secondarg\csname#1\endcsname}\fi}

\def\pagerefs#1{\expandafter\ifx\csname#1\endcsname\relax
  {\nodefmsg{#1}\bf`#1'}\else
  \expandafter\firstarg\csname#1\endcsname%
  s~\htmllocref{#1}{\expandafter\secondarg\csname#1\endcsname}\fi}

\def\pagerefn#1{\expandafter\ifx\csname#1\endcsname\relax
  {\nodefmsg{#1}\bf`#1'}\else
  \htmllocref{#1}{\expandafter\secondarg\csname#1\endcsname}\fi}






\newif\ifhypers  
\hypersfalse


\edef\freehash{\catcode`\noexpand\#=\the\catcode`\#}%
\catcode`\#=12
\freehash
\let\freehash=\relax
\ifhypers\fi
\def\puthtml#1{\ifhypers\fi}
\def\htmlanchor#1#2{\puthtml{<a name="#1">}#2\puthtml{</a>}}
\def\@pdfm@mark#1{}
\def\setlink#1{\colored{\linkcolor}{#1}}%
\def\setlink#1{\ifdraft\Purple{#1}\else{#1}\fi}%

%

%
\def\htmllocref#1#2{\ifhypers\leavevmode\fi\setlink{#2}\ifhypers\fi\relax}%
%
%
%
%
\def\Acrobatmenu#1#2{%
  \@pdfm@mark{%
    bann <<
      /Type /Annot
      /Subtype /Link
      /A <<
        /S /Named
        /N /#1
      >>
      /Border [\@pdfborder]
      /C [\@menubordercolor]
    >>%
   }%
  \Hy@colorlink{\@menucolor}#2\Hy@endcolorlink
  \@pdfm@mark{eann}%
}
\def\@pdfborder{0 0 1}
\def\@menubordercolor{1 0 0}
\def\@menucolor{red}

\def\ifempty#1#2\endB{\ifx#1\endA}
\def\makeref#1#2#3{\ifempty#1\endA\endB\else\forward{#1}{#2}{#3}\fi\unskip}


\def\crosshairs#1#2{
  \dimen1=.002\drawx
  \dimen2=.002\drawy
  \ifdim\dimen1<\dimen2\dimen3\dimen1\else\dimen3\dimen2\fi
  \setbox1=\vclap{\vline{2\dimen3}}
  \setbox2=\clap{\hline{2\dimen3}}
  \setbox3=\hstutter{\kern\dimen1\box1}{4}
  \setbox4=\vstutter{\kern\dimen2\box2}{4}
  \setbox1=\vclap{\vline{4\dimen3}}
  \setbox2=\clap{\hline{4\dimen3}}
  \setbox5=\clap{\copy1\hstutter{\box3\kern\dimen1\box1}{6}}
  \setbox6=\vclap{\copy2\vstutter{\box4\kern\dimen2\box2}{6}}
  \setbox1=\vbox{\offinterlineskip\box5\box6}
  \smash{\vbox to #2{\hbox to #1{\hss\box1}\vss}}}

\def\boxgrid#1{\rlap{\vbox{\offinterlineskip
  \setbox0=\hline{\wd#1}
  \setbox1=\vline{\ht#1}
  \smash{\vbox to \ht#1{\offinterlineskip\copy0\vfil\box0}}
  \smash{\vbox{\hbox to \wd#1{\copy1\hfil\box1}}}}}}

\def\drawgrid#1#2{\vbox{\offinterlineskip
  \dimen0=\drawx
  \dimen1=\drawy
  \divide\dimen0 by 10
  \divide\dimen1 by 10
  \setbox0=\hline\drawx
  \setbox1=\vline\drawy
  \smash{\vbox{\offinterlineskip
    \copy0\vstutter{\kern\dimen1\box0}{10}}}
  \smash{\hbox{\copy1\hstutter{\kern\dimen0\box1}{10}}}}}

\long\def\boxit#1#2#3{\hbox{\vrule
  \vtop{%
    \vbox{\hsize=#2\hrule\kern#1%
      \hbox{\kern#1#3\kern#1}}%
    \kern#1\hrule}%
    \vrule}}
\long\def\boxitv#1#2#3{\boxit{#1}{#2}{{\hbox to #2{\vrule%
  \vbox{\hsize=#2\hrule\kern#1%
      \hbox{\kern#1#3\kern#1}}%
    \kern#1\hrule}%
    \vrule}}}

\newdimen\boxingdimen
\long\def\boxing#1{\boxingdimen=\hsize\advance\boxingdimen by -2ex%
\vskip0.5cm\boxit{1ex}{\boxingdimen}{\vbox{#1}}\vskip0.5cm}

\newdimen\boxrulethickness \boxrulethickness=.4pt
\newdimen\boxedhsize
\newbox\textbox
\newdimen\originalbaseline
\newdimen\hborderspace\newdimen\vborderspace
\hborderspace=3pt \vborderspace=3pt

\def\preparerulebox#1#2{\setbox\textbox=\hbox{#2}%
   \originalbaseline=\vborderspace
   \advance\originalbaseline\boxrulethickness
   \advance\originalbaseline\dp\textbox
   \def\Borderbox{\vbox{\hrule height\boxrulethickness
     \hbox{\vrule width\boxrulethickness\hskip\hborderspace
     \vbox{\vskip\vborderspace\relax#1{#2}\vskip\vborderspace}%
     \hskip\hborderspace\vrule width\boxrulethickness}%
     \hrule height\boxrulethickness}}}

\def\hrulebox#1{\hbox{\preparerulebox{\hbox}{#1}%
   \lower\originalbaseline\Borderbox}}
\def\vrulebox#1#2{\vbox{\preparerulebox{\vbox}{\hsize#1#2}\Borderbox}}
\def\parrulebox#1\par{\boxedhsize=\hsize\advance\boxedhsize
   -2\boxrulethickness \advance\boxedhsize-2\hborderspace
   \vskip3\parskip\vrulebox{\boxedhsize}{#1}\vskip\parskip\par}
\def\sparrulebox#1\par{\vskip1truecm\boxedhsize=\hsize\advance\boxedhsize
   -2\boxrulethickness \advance\boxedhsize-2\hborderspace
   \vskip3\parskip\vrulebox{\boxedhsize}{#1}\vskip.5truecm\par}



%

\newcount\sectno \sectno=0
\newcount\appno \appno=64 
\newcount\thmno \thmno=0
\newcount\randomct \newcount\rrandomct
\newif\ifsect \secttrue
\newif\ifapp \appfalse
\newif\ifnumct \numcttrue
\newif\ifchap \chapfalse
\newif\ifdraft \draftfalse
\newif\ifcolor\colorfalse


\def\widow#1{\vskip 0pt plus#1\vsize\goodbreak\vskip 0pt plus-#1\vsize}


\def\stdskip{\vskip 9pt plus3pt minus 3pt}
\def\stdbreak{\par\removelastskip\penalty-100\stdskip}
\def\halfbreak{\vskip 0.6ex\penalty-100}

\def\proof{\stdbreak\noindent%
  \ifdraft\let\proofcolor=\Purple\else\let\proofcolor=\Black\fi%
  \proofcolor{{\sl Proof. }}}

\def\proofone{\stdbreak\noindent%
  \ifdraft\let\proofcolor=\Purple\else\let\proofcolor=\Black\fi%
  \proofcolor{{\sl Proof \#1:\ \ }}}
\def\prooftwo{\stdbreak\noindent%
  \ifdraft\let\proofcolor=\Purple\else\let\proofcolor=\Black\fi%
  \proofcolor{{\sl Proof \#2:\ \ }}}
\def\proofthree{\stdbreak\noindent%
  \ifdraft\let\proofcolor=\Purple\else\let\proofcolor=\Black\fi%
  \proofcolor{{\sl Proof \#3:\ \ }}}
\def\proofof#1{\stdbreak\noindent%
  \ifdraft\let\proofcolor=\Purple\else\let\proofcolor=\Black\fi%
  \proofcolor{{\sl Proof of #1:\ }}}
\def\claim{\stdbreak\noindent%
  \ifdraft\let\proofcolor=\Purple\else\let\proofcolor=\Black\fi%
  \proofcolor{{\sl Claim:\ }}}
\def\proofclaim{\stdbreak\noindent%
  \ifdraft\let\proofcolor=\Purple\else\let\proofcolor=\Black\fi%
  \proofcolor{{\sl Proof of the claim:\ }}}

\newif\ifnumberbibl \numberbibltrue
\newdimen\labelwidth

\def\references{\bgroup
  \edef\numtoks{}%
  \global\thmno=0
  \setbox1=\hbox{[999]} 
  \labelwidth=\wd1 \advance\labelwidth by 2.5em
  \ifnumberbibl\advance\labelwidth by -2em\fi
  \parindent=\labelwidth \advance\parindent by 0.5em 
  \removelastskip
  \widow{0.1}
  \vskip 24pt plus 6pt minus 6 pt
  \frenchspacing
  \immediate\write\isauxout{\noexpand\forward{Bibliography}{}{\the\pageno}}%
  \immediate\write\iscontout{\noexpand\contnosectlist{}{Bibliography}{\the\pageno}}%
  \ifdraft\let\refcolor=\Maroon\else\let\refcolor=\Black\fi%
  \leftline{\sectionfont\refcolor{References}}
  \ifhypers%
     \global\thmno=0\relax%
     \ifsect%
	\global\advance\sectno by 1%
	\edef\numtoks{\number\sectno}%
  	\hbox{}%
     \else%
	\edef\numtoks{}%
  	\hbox{}%
     \fi%
  \fi%
  \nobreak\stdskip\noindent}%
\def\endreferences{\nonfrenchspacing\egroup}

\def\referencesn{\bgroup
  \edef\numtoks{}%
  \global\thmno=0
  \setbox1=\hbox{[999]} 
  \labelwidth=\wd1 \advance\labelwidth by 2.5em
  \ifnumberbibl\advance\labelwidth by -2em\fi
  \parindent=\labelwidth \advance\parindent by 0.5em 
  \removelastskip
  \widow{0.1}
  \vskip 24pt plus 6pt minus 6 pt
  \frenchspacing
  \immediate\write\isauxout{\noexpand\forward{Bibliography}{}{\the\pageno}}%
  \immediate\write\iscontout{\noexpand\contnosectlist{}{Bibliography}{\the\pageno}}%
  \ifdraft\let\refcolor=\Maroon\else\let\refcolor=\Black\fi%
  \ifhypers%
     \global\thmno=0\relax%
     \ifsect%
	\global\advance\sectno by 1%
	\edef\numtoks{\number\sectno}%
  	\hbox{}%
     \else%
	\edef\numtoks{}%
  	\hbox{}%
     \fi%
  \fi%
  \nobreak\stdskip\noindent}%

\def\bitem#1{\global\advance\thmno by 1%
  \ifdraft\let\itemcolor=\Purple\else\let\itemcolor=\Black\fi%
  \outer\expandafter\def\csname#1\endcsname{\the\thmno}%
  \edef\numtoks{\number\thmno}%
  \ifhypers\htmlanchor{#1}{\makeref{\noexpand#1}{REF}{\numtoks}}\fi%
  \ifnumberbibl
    \immediate\write\isauxout{\noexpand\forward{\noexpand#1}{}{\the\thmno}}
    \item{\hbox to \labelwidth{\itemcolor{\hfil\the\thmno.\ \ }}}
  \else
    \immediate\write\isauxout{\noexpand\expandafter\noexpand\gdef\noexpand\csname #1\noexpand\endcsname{#1}}
    \item{\hbox to \labelwidth{\itemcolor{\hfil#1\ \ }}}
  \fi}

\newcount\chapno
\newcount\chappage
\chapno=0

\newtoks\rightpagehead
\footline={\bgroup{\ifdraft\ifnum\pageno>5\hfill\hbox{\headerfont\today}\fi\fi}\egroup}

\outer\def\section#1{%
  \widow{0.1}%
  \global\subsectno=0%
  \appfalse\secttrue%
  \removelastskip%
  \global\advance\sectno by 1%
  \global\thmno=0\relax%
  \global\randomct=0\relax%
  \edef\numtoks{\ifchap\number\chapno.\fi\number\sectno}%
  \edef\secttitl{#1}%
  \edef\secttitle{Section \numtoks: \secttitl}%
  \rightpagehead={\Black{\headerfont\hfill\secttitle\kern3em\the\pageno}}%
  \vskip 24pt plus 6pt minus 6 pt%
  \ifdraft\let\sectcolor=\OliveGreen\else\let\sectcolor=\Black\fi%
  \message{#1}%
  \ifhypers\hbox{}\fi%
    \futurelet\testchar\maybeoptionsection}

\def\maybeoptionsection{\ifx[\testchar\let\next\optionsection%
	\else\let\next\nooptionsection\fi\next}

\def\optionsection[#1]{%
  \immediate\write\isauxout{\noexpand\forward{\noexpand#1}{Section}{\numtoks}}%
  {\noindent{\draftlabel{#1}\sectionfont\sectcolor{\numtoks}\quad \sectcolor{\secttitl}}}%
  \immediate\write\iscontout{\noexpand\contlist{\noexpand#1}{\secttitl}{\the\pageno}}%
  \htmlanchor{#1}{\makeref{#1}{Section}{\numtoks}}%
  \nobreak\vskip 4ex}

\def\nooptionsection{%
  {\noindent{\sectionfont\sectcolor{\numtoks}\quad \sectcolor{\secttitl}}}%
  \write\iscontout{\noexpand\contlist{0}{\secttitle}{\the\pageno}}%
  \nobreak\vskip 4ex}

\newcount\subsectno
\outer\def\subsection#1{%
  \advance\subsectno by 1%
  \global\randomct=0\relax%
  \vskip 10pt plus 6pt minus 6 pt%
  \widow{.02}%
  \message{#1}%
  \noindent \sectcolor{$\underline{\hbox{\numtoks.\number\subsectno\quad #1}}$}%
  \nobreak\vskip 1ex}


\newif\ifnewgroup\newgroupfalse

%
\def\proclamationsing#1#2#3#4{
  \outer\expandafter\def\csname#1\endcsname{%
  \global\newgroupfalse%
  \ifnum#3<5\global\newgrouptrue\fi%
  \ifnum#3<1\global\newgroupfalse\fi%
  \ifdraft\let\proclaimcolor=\Fuchsia\else\let\proclaimcolor=\Black\fi%
  \gdef\Prnm{#4}%
  \global\advance\thmno by 1%
  \global\randomct=0%
  \ifcase#3
	\stdbreak \or 
	\stdbreak \or 
	\stdbreak \or 
	\stdbreak \or 
	\halfbreak \or 
	\halfbreak \or 
	\halfbreak \or 
	\halfbreak \or 
	\halfbreak \or 
	\halfbreak \or 
	\halfbreak \or 
	\else \stdbreak \fi%
  \edef\proctoks{\ifchap\the\chapno.\fi\ifsect\the\sectno.\fi\the\thmno\ifnum#3=2'\fi\ifnum#3=3''\fi}%
  \widow{0.05}%
  \ifnumct\noindent{\proclaimcolor{%
	\ifnum#3=8*\fi%
	\ifnum#3=9*\fi%
	\ifnum#3=7\hbox to 1ex{$\dag$}\fi%
	\ifnum#3=10\hbox to 0.5ex{$\dag$}\fi%
	\ifnum#3=-1{\bf Important}\ \lowercase\fi%
	\ifnum#3=6\else\ifnum#3=9\else\ifnum#3=10\else{\bf #2}\fi\fi\fi%
	\ifnum#3=6\else\ifnum#3=9\else\ \fi\fi%
	\bf \proctoks}\enspace}%
  \else\noindent{\proclaimcolor{\bf #2}\enspace}%
  \fi%
  \futurelet\testchar\maybeoptionproclaim}}

\def\proclamation#1#2#3{
  \outer\expandafter\def\csname#1\endcsname{%
  \global\newgroupfalse%
  \ifnum#3<5\global\newgrouptrue\fi%
  \ifnum#3<1\global\newgroupfalse\fi%
  \ifdraft\let\proclaimcolor=\Fuchsia\else\let\proclaimcolor=\Black\fi%
  \gdef\Prnm{#2}%
  \global\advance\thmno by 1%
  \global\randomct=0%
  \ifcase#3
	\stdbreak \or 
	\stdbreak \or 
	\stdbreak \or 
	\stdbreak \or 
	\halfbreak \or 
	\halfbreak \or 
	\halfbreak \or 
	\halfbreak \or 
	\halfbreak \or 
	\halfbreak \or 
	\halfbreak \or 
	\else \stdbreak \fi%
  \edef\proctoks{\ifchap\the\chapno.\fi\ifsect\the\sectno.\fi\the\thmno\ifnum#3=2'\fi\ifnum#3=3''\fi}%
  \widow{0.05}%
  \ifnumct\noindent{\proclaimcolor{%
	\ifnum#3=8*\fi%
	\ifnum#3=9*\fi%
	\ifnum#3=7\hbox to 1ex{$\dag$}\fi%
	\ifnum#3=10\hbox to 0.5ex{$\dag$}\fi%
	\ifnum#3=-1{\bf Important}\ \lowercase\fi%
	\ifnum#3=6\else\ifnum#3=9\else\ifnum#3=10\else{\bf #2}\fi\fi\fi%
	\ifnum#3=6\else\ifnum#3=9\else\ \fi\fi%
	\bf \proctoks}\enspace}%
  \else\noindent{\proclaimcolor{\bf #2}\enspace}%
  \fi%
  \futurelet\testchar\maybeoptionproclaim}}

\def\maybeoptionproclaim{\ifx[\testchar\let\next\optionproclaim%
	\else\let\next\nooptionproclaim\fi\next}

\def\optionproclaim[#1]{%
  \bgroup\htmlanchor{#1}{\makeref{\noexpand#1}{\Prnm}{\proctoks}}\egroup%
  \immediate\write\isauxout{\noexpand\forward{\noexpand#1}{\Prnm}{\proctoks}}%
  \draftlabel{#1}%
  \ifnewgroup\bgroup\sl\fi}
\def\nooptionproclaim{\ifnewgroup\bgroup\sl\fi}
\def\endb{\par\stdbreak\egroup\newgroupfalse\randomct=0}

\def\pagelabel#1{%
   \ \unskip%
   \write\isauxout{\noexpand\forward{\noexpand#1}{page}{\the\pageno}}%
   \bgroup\htmlanchor{#1}{\makeref{\noexpand#1}{page}{\the\pageno}}\egroup%
   \draftlabel{page\ #1}%
 }

\def\plot[#1]{
    \ifdraft\let\proclaimcolor=\Fuchsia\else\let\proclaimcolor=\Black\fi%
    \gdef\Prnm{Plot}%
    \global\advance\thmno by 1%
    \ifnumct{\bf \proclaimcolor{Plot\ \proctoks}}%
    \else\noindent{\bf \proclaimcolor{Plot}}%
    \fi%
  \htmlanchor{#1}{\makeref{\noexpand#1}{\Prnm}{\proctoks}}%
  \immediate\write\isauxout{\noexpand\forward{\noexpand#1}{\Prnm}{\proctoks}}%
  \draftlabel{#1}%
  }

\newcount\captionct \captionct=0
\outer\def\caption#1{{\global\advance\captionct by 1%
  \relax%
  \gdef\Irnm{Figure}%
  \edef\numtoks{\ifchap\the\chapno.\fi\ifsect\the\sectno.\fi\the\captionct}%
  \vskip0.5cm\vtop{\noindent{\bf Figure \numtoks.} #1}\vskip0.5cm%
  \futurelet\testchar\maybeoptionproclaim}}
\outer\def\caption#1{{\global\advance\thmno by 1%
  \relax%
  \gdef\Irnm{Figure}%
  \edef\proctoks{\ifchap\the\chapno.\fi\ifsect\the\sectno.\fi\the\thmno}%
  \vskip0.2cm\vtop{\narrower\narrower\noindent{\bf Figure \proctoks.} #1}\vskip0.5cm%
  \futurelet\testchar\maybeoptionproclaim}}

\def\eqlabel#1{
    \ifdraft\let\labelcolor=\Red\else\let\labelcolor=\Black\fi%
    \global\advance\thmno by 1%
    \edef\proctoks{\ifchap\number\chapno.\fi\ifsect\number\sectno.\fi\number\thmno}%
    \outer\expandafter\def\csname#1\endcsname{\proctoks}%
    \htmlanchor{#1}{\makeref{\noexpand#1}{}{\proctoks}}%
    \immediate\write\isauxout{\noexpand\forward{\noexpand#1}{Equation}{(\proctoks)}}%
    \quad\qquad\hfill\eqno\hbox{\hfill\rm\noexpand\labelcolor{(\proctoks)}\draftcmt{\labelcolor{#1}}}%
}
\def\iqlabel#1{
    \ifdraft\let\labelcolor=\Red\else\let\labelcolor=\Black\fi%
    \global\advance\thmno by 1%
    \edef\proctoks{\ifchap\number\chapno.\fi\ifsect\number\sectno.\fi\number\thmno}%
    \outer\expandafter\def\csname#1\endcsname{\proctoks}%
    \htmlanchor{#1}{\makeref{\noexpand#1}{}{\proctoks}}%
    \immediate\write\isauxout{\noexpand\forward{\noexpand#1}{Inequality}{(\proctoks)}}%
    \quad\qquad\hfill\eqno\hbox{\hfill\rm\noexpand\labelcolor{(\proctoks)}\draftcmt{\labelcolor{#1}}}%
}

\def\eqalabel#1{
    \ifdraft\global\let\labelcolor=\Red\else\global\let\labelcolor=\Black\fi%
    \global\advance\thmno by 1%
    \gdef\proctoks{\ifchap\number\chapno.\fi\ifsect\number\sectno.\fi\number\thmno}%
    \outer\expandafter\def\csname#1\endcsname{\proctoks}%
    \htmlanchor{#1}{\makeref{\noexpand#1}{}{\proctoks}}%
    \immediate\write\isauxout{\noexpand\forward{\noexpand#1}{equation}{(\proctoks)}}%
    \quad\qquad&\hfill\hbox{\hfill\rm\noexpand\labelcolor{(\proctoks)}\draftcmt{\labelcolor{#1}}}%
}

\def\label#1{\optionproclaim[#1]}
\def\label#1{%
    \ifdraft\let\labelcolor=\Red\else\let\labelcolor=\Black\fi%
    \edef\proctoks{\ifchap\the\chapno.\fi\ifsect\the\sectno.\fi\the\thmno}%
    \outer\expandafter\def\csname#1\endcsname{\proctoks}%
    \bgroup\htmlanchor{#1}{\makeref{\noexpand#1}{}{\proctoks}}\egroup%
    \immediate\write\isauxout{\noexpand\forward{\noexpand#1}{Figure}{\proctoks}}%
}

\proclamation{definition}{Definition}{1}
\proclamation{defin}{Definition}{1}
\proclamation{defins}{Definitions}{1}
\proclamation{lemma}{Lemma}{1}
\proclamation{lemmap}{Lemma}{2}
\proclamation{lemmapp}{Lemma}{3}
\proclamation{proposition}{Proposition}{1}
\proclamation{prop}{Proposition}{1}
\proclamation{theorem}{Theorem}{1}
\proclamation{thm}{Theorem}{1}
\proclamation{corollary}{Corollary}{1}
\proclamation{cor}{Corollary}{1}
\proclamation{conjecture}{Conjecture}{1}
\proclamation{proc}{Procedure}{1}
\proclamation{assumption}{Assumption}{1}
\proclamation{axiom}{Axiom}{1}
\proclamation{example}{Example}{0}
\proclamation{examples}{Examples}{0}
\proclamation{examplesnotation}{Examples and notation}{0}
\proclamation{exnoex}{Examples and non-examples}{0}
\proclamation{remark}{Remark}{0}
\proclamation{remarks}{Remarks}{0}
\proclamation{remarksexamples}{Remarks and examples}{0}
\proclamationsing{facts}{Facts}{0}{Fact}
\proclamation{fact}{Fact}{0}
\proclamation{irrelevant}{Irrelevant}{0}
\proclamation{question}{Question}{0}
\proclamation{construction}{Construction}{0}
\proclamation{algorithm}{Algorithm}{0}
\proclamation{problem}{Problem}{0}
\proclamation{table}{Table}{0}
\proclamation{notation}{Notation}{0}
\proclamation{figure}{Figure}{0}
\proclamation{impexample}{Example}{-1}
\proclamation{exercise}{Exercise}{0}
\proclamation{importantexercise}{IMPORTANT EXERCISE}{0}
\proclamation{exercisedag}{Exercise}{7} 
\proclamation{stexercise}{Exercise}{8} 
\proclamation{excise}{Exercise}{6} 
\proclamation{stexcise}{Exercise}{9} 
\proclamation{excisedag}{Exercise}{10} 

\def\doubleindent{\advance \leftskip 2\parindent\advance \rightskip \parindent}
\def\singleindent{\advance \leftskip \parindent\advance \rightskip \parindent}
\def\doubleleftindent{\advance \leftskip 2\parindent}


\def\hpad#1#2#3{\hbox{\kern #1\hbox{#3}\kern #2}}
\def\vpad#1#2#3{\setbox0=\hbox{#3}\vbox{\kern #1\box0\kern #2}}



\def\stack#1#2#3{\vbox{\offinterlineskip
  \setbox2=\hbox{#2}
  \setbox3=\hbox{#3}
  \dimen0=\ifdim\wd2>\wd3\wd2\else\wd3\fi
  \hbox to \dimen0{\hss\box2\hss}
  \kern #1
  \hbox to \dimen0{\hss\box3\hss}}}


\def\hexp#1{%
  \setbox0=\hbox{${}^{#1}$}%
  \hbox to .5\wd0{\box0\hss}}

\def\hsub#1{%
  \setbox0=\hbox{${}_{#1}$}%
  \hbox to .5\wd0{\box0\hss}}


\def\mOth{\mathsurround=0pt}
\newdimen\pOrenwd \setbox0=\hbox{\tenex B} \pOrenwd=\wd0
\def\leftextramatrix#1{\begingroup \mOth
  \setbox0=\vbox{\def\cr{\crcr\noalign{\kern2pt\global\let\cr=\endline}}
    \ialign{$##$\hfil\kern2pt\kern\pOrenwd&\thinspace\hfil$##$\hfil
      &&\quad\hfil$##$\hfil\crcr
      \omit\strut\hfil\crcr\noalign{\kern-\baselineskip}
      #1\crcr\omit\strut\cr}}
  \setbox2=\vbox{\unvcopy0 \global\setbox1=\lastbox}
  \setbox2=\hbox{\unhbox1 \unskip \global\setbox1=\lastbox}
  \setbox2=\hbox{$\kern\wd1\kern-\pOrenwd \left[ \kern-\wd1
    \vcenter{\unvbox0 \kern-\baselineskip} \,\right]$}
  \null\;\vbox{\box2}\endgroup}


\def\COMMENT#1\par{\bigskip\hrule\smallskip#1\smallskip\hrule\bigskip}

{\catcode`\^^M=12 \endlinechar=-1 %
 \gdef\xcomment#1^^M{\def\test{#1}
   \ifx\test\endcomment \let\next=\endgroup
   \else \let\next=\xcomment \fi
   \next}
}
\def\dospecials{\do\ \do\\\do\{\do\}\do\$\do\&%
  \do\#\do\^\do\^^K\do\_\do\^^A\do\%\do\~}
\def\makeinnocent#1{\catcode`#1=12}
\def\comment{\begingroup
  \let\do=\makeinnocent \dospecials
  \endlinechar`\^^M \catcode`\^^M=12 \xcomment}
{\escapechar=-1
 \xdef\endcomment{\string\\endcomment}
}



\font\foott=cmtt9
\def\draftcmt#1{\ifdraft\llap{\smash{\raise2ex\hbox{{#1}}}}\fi}
\def\draftind#1{\llap{\smash{\raise-1.3ex\hbox{{#1}}}}}
\def\draftlabel#1{\ifdraft\draftcmt{\footfonts\Green{{\string#1}\ }\unskip\fi}}
\def\draftlabel#1{%
  \ifx#1\empty\else
  \ifdraft\draftcmt{\foott\Green{\string #1}\ \ }\unskip\fi\fi}

\newbox\strikebox
\def\strike#1{\setbox\strikebox=\hbox{#1}%
  \rlap{\raise0.4ex\hbox to \wd\strikebox{\leaders\hrule height 0.8pt depth 0pt\hfill}}%
  #1}
\def\strikem#1{\setbox\strikebox=\hbox{$#1$}%
  \rlap{\raise0.4ex\hbox to \wd\strikebox{\leaders\hrule height 0.8pt depth 0pt\hfill}}%
  \hbox{$#1$}}

\newcount\hours
\newcount\minutes
\hours\time \divide\hours 60
\minutes-\hours \multiply\minutes 60\advance\minutes\time
\edef\timehhmm{\ifnum\hours<10 0\fi\the\hours
:\ifnum\minutes<10 0\fi\the\minutes}

\footline={{\rm\hfill\folio\hfill}}

\def\cite{\futurelet\testchar\maybeoptioncite}
\def\maybeoptioncite{\ifx[\testchar \let\next\optioncite
	\else \let\next\nooptioncite \fi \next}

\def\nooptioncite#1{\expandafter\ifx\csname#1\endcsname\relax
  {\nodefmsg{#1}\bf`#1'}\else
  \htmllocref{#1}{[\csname#1\endcsname]}{}\fi}

\def\optioncite[#1]#2{\ifx\csname#1\endcsname\relax\else
  \htmllocref{#2}{[\csname#2\endcsname, #1]}{}\fi}


\def\strutdepth{\dp\strutbox}
\def\strutdepth{4.3ex}
\def\margin#1{{\strut\vadjust{\kern-\strutdepth%
        \vtop to\strutdepth{\baselineskip\strutdepth%
        \llap{\vbox{\hsize=3em\noindent #1}\hskip1em}\null}}}}
\def\rmargin#1{\strut\vadjust{\kern-\strutdepth%
        \vtop to\strutdepth{\baselineskip\strutdepth%
        \vss\rlap{\ \hskip\hsize\ \vtop{\hsize=4em\noindent #1}}\null}}}
\def\warning{\vskip0ex\noindent\llap{\lower4pt\hbox to 1.0cm{\psline(0,.3)(.27,.3)(0,.09)(.03,0)(.27,0)(.3,.03)(.27,.06)}\hfill\ }\indent}
\def\warning{\margin{\lower4pt\hbox to 0.6cm{\psline(0,.3)(.27,.3)(0,.09)(.03,0)(.27,0)(.3,.03)(.27,.06)}\hfill\ }}
\def\warning{\margin{\hbox to 0.6cm{\psline(0,.3)(.27,.3)(0,.09)(.03,0)(.27,0)(.3,.03)(.27,.06)}\hfill\ }}
\def\warning{\noindent{\bf Warning: }\margin{\ \ \ \ \lower-4pt\hbox to 0.6cm{\psline[linewidth=2pt](0,0)(.5,0)(0.53,-0.02)(0.55,-0.05)(0.53,-0.1)(0.03,-0.35)(0.01,-0.4)(0.03,-0.43)(0.06,-0.45)(0.56,-0.45)}\hfill}}
\def\warningdown#1{\noindent{\bf Warning: }\bgroup\def\strutdepth{#1}\margin{\ \ \
\ \hbox to 0.6cm{\psline[linewidth=2pt](0,0)(.5,0)(0.53,-0.02)(0.55,-0.05)(0.53,-0.1)(0.03,-0.35)(0.01,-0.4)(0.03,-0.43)(0.06,-0.45)(0.56,-0.45)}\hfill}\egroup}

\def\quote{\bgroup\footfont\parindent=0pt%
	\baselineskip=8pt\hangindent=15em\hangafter=0}
\def\endquote{\egroup}

\def\bmargin#1{{\strut\vadjust{\kern-\strutdepth%
        \vtop to\strutdepth{\baselineskip\strutdepth%
        \llap{\vbox{\hsize=3em\noindent #1}\hskip2em}\null}}}}

\newdimen\hngind \hngind=\parindent \advance\hngind by 3.9em\relax
\def\nonuritem{{\vskip0pt\hangafter=0\global\hangindent\hngind%
        {\vphantom{(\the\randomct)\hskip.5em}}}}
\def\nonuritem{{\parindent=0pt\hngind=\parindent \advance\hngind by 1.9em%
        \vskip0pt\hangafter=0\global\hangindent\hngind%
        \vphantom{{(\the\randomct)\hskip.5em}}}}
\def\nonuritem{{\idoitem{}}}
\def\pitem{{\global\advance\randomct by 1%
	\vskip1pt\hangafter=0\global\hangindent\hngind%
	{{\the\randomct.\hskip.5em}}}}

%
\def\idoitem#1{{\parindent=0pt\hngind=\parindent \advance\hngind by 3.7em%
        \vskip0pt\hangafter=0\global\hangindent\hngind%
        {$\hbox to3.7em{\hfill#1\hskip.5em}$}}}%
\def\romanitem{\global\advance\randomct by 1%
    \edef\itoks{(\romannumeral\randomct)}%
    \edef\Irnm{}%
    \idoitem{\itoks}%
    \futurelet\testchar\maybeoptioniitem} 
\def\sitem{\global\advance\randomct by 1%
    \edef\itoks{(\the\randomct)}%
    \edef\Irnm{Step}%
    \idoitem{\itoks}%
    \futurelet\testchar\maybeoptioniitem} 
\def\Rsitem{\global\advance\rrandomct by 1%
    \edef\itoks{(\uppercase\expandafter{\romannumeral\rrandomct)}}%
    \edef\Irnm{Step}%
    \idoitem{\itoks}%
    \futurelet\testchar\maybeoptioniitem} 
\def\citem{\global\advance\randomct by 1%
    \edef\itoks{(\the\randomct)}%
    \edef\Irnm{Condition}%
    \idoitem{\itoks}%
    \futurelet\testchar\maybeoptionitem} 
\def\iitem{\global\advance\randomct by 1%
    \edef\itoks{(\the\randomct)}%
    \edef\Irnm{}%
    \idoitem{\itoks}%
    \futurelet\testchar\maybeoptioniitem} 
\def\rpitem{\global\advance\randomct by 1%
    \edef\itoks{(\the\randomct')}%
    \edef\Irnm{}%
    \idoitem{\itoks}%
    \futurelet\testchar\maybeoptioniitem} 

\def\maybeoptioniitem{\ifx[\testchar\let\next\optioniitem%
	\else\let\next\nooptioniitem\fi\next}
\def\optioniitem[#1]{%
  \bgroup\htmlanchor{#1}{\makeref{\noexpand#1}{\Irnm}{\itoks}}\egroup%
  \immediate\write\isauxout{\noexpand\forward{\noexpand#1}{\Irnm}{\itoks}}%
  \draftlabel{#1}}
\def\nooptioniitem{}
\def\bulletpoint{{\parindent=0pt\hngind=\parindent \advance\hngind by 1.2em%
        \vskip0pt\hangafter=0\global\hangindent\hngind%
        {$\hbox to1.2em{\hfill$\bullet$\hskip.5em}$}}}%
\def\bulletpoint{{\hngind=\parindent \advance\parindent by -1.2em%
        \vskip0pt\hangafter=0\global\hangindent\hngind%
        {$\hbox to1.2em{\hfill$\bullet$\hskip.5em}$}}}%
\def\nopoint{{\parindent=0pt\hngind=1cm \advance\hngind by 1.2em%
        \vskip0pt\hangafter=0\global\hangindent\hngind%
        {$\hbox to1.2em{\hfill\ \hskip.5em}$}}}%
\def\nopoint{{\hngind=\parindent \advance\parindent by -1.2em%
        \vskip0pt\hangafter=0\global\hangindent\hngind%
        {$\hbox to1.2em{\hfill\ \hskip.5em}$}}}%
\def\bbulletpoint{{\bgroup\parskip=0.1cm\parindent=1.5cm\hngind=\parindent%
	\advance\hngind by 1.5em%
        \vskip0pt\hangafter=0\global\hangindent\hngind%
        {$\hbox to1.5em{\hfill$\bullet$\hskip.5em}$}\egroup}}%
\def\bnobulletpoint{{\bgroup\parskip=0.1cm\parindent=1.5cm\hngind=\parindent%
	\advance\hngind by 1.5em%
        \vskip0pt\hangafter=0\global\hangindent\hngind%
        {$\hbox to1.5em{\hfill\hskip.5em}$}\egroup}}%

\let\ritem=\iitem

\def\rlabel#1{{%
  \bgroup\htmlanchor{#1}{\makeref{\noexpand#1}{item}{\the\randomct}}\egroup%
  \immediate\write\isauxout{\noexpand\forward{\noexpand#1}{item}{\the\randomct}}%
  \draftlabel{item\ #1}}}

\def\idoitemitem#1{{\hngind=\parindent \advance\hngind by 4.3em%
        \vskip-\parskip\hangafter=0\global\hangindent\hngind%
        {$\hbox to3.9em{\hfill#1\ }$}}}
\def\idoitemitemm#1{{\hngind=\parindent \advance\hngind by 7.3em%
        \vskip-\parskip\hangafter=0\global\hangindent\hngind%
        {$\hbox to6.9em{\hfill#1\ }$}}}
\let\riitem=\idoitemitem

\def\thmparens#1{\ignorespaces{\rm(#1)}}

\newdimen\thmskip \thmskip=0.7ex
	
\def\thmalready#1{\removelastskip\vskip\thmskip%
	\global\randomct=0%
	\noindent{\bf #1:}\ \ %
	\bgroup \it%
	\abovedisplayskip=4pt\belowdisplayskip=3pt%
	\parskip=0pt%
	}


\newdimen\exerciseindent \exerciseindent=4.5em
\newdimen\exerciseitemindent \exerciseitemindent=1.9em
\newdimen\exitindwk \exitindwk=3.7em
\def\exitem{\global\advance\randomct by 1%
	\vskip0pt\indent\hangindent\exitindwk \hskip\exerciseitemindent%
	\llap{\romannumeral\randomct)\hskip0.5em}\ \hskip-0.5em}
\def\stexitem{\global\advance\randomct by 1%
	\par\indent\hangindent\exitindwk \hskip\exerciseitemindent%
	\llap{(\romannumeral\randomct)\hskip-.1ex*\enspace}\ignorespaces}


\def\exercises{\global\thmno=0\vskip0.9cm%
	\noindent{\bf Exercises for Section~\ifchap\the\chapno.\fi\the\sectno}}
\def\solutions{\global\thmno=0\bigskip\widow{0.05}%
	\noindent{\bf Solutions for Section~\ifchap\the\chapno.\fi\the\sectno}}
\def\solution#1{\removelastskip\vskip 5pt%
	\global\randomct=0%
	\parindent=\exerciseindent%
	\par\hangindent\exerciseindent\indent\llap{\hbox to \exerciseindent{%
	  {\bf \unskip\ref{#1}}}}}
\def\solution#1{\removelastskip\vskip 5pt\par\noindent {\bf \refn{#1}:}}


\newdimen\wddim
\def\edinsert#1{\setbox0=\hbox{#1}\wddim=\wd0\divide\wddim by 2%
\smash{\hskip-\wddim\raise15pt\hbox to 0pt{\Red{$\underbrace{\box0}$}}\hskip\wddim}}

\def\midline#1{\raise 0.20em\hbox to #1{\vrule depth0pt height 0.405pt width #1}}

\def\buildrelu#1\over#2{\mathrel{\mathop{\kern0pt #1}\limits_{#2}}}



\def\bmatrix#1#2{{\left(\vcenter{\halign
  {&\kern#1\hfil$##\mathstrut$\kern#1\cr#2}}\right)}}

\def\rightarrowmat#1#2#3{
  \setbox1=\hbox{\small\kern#2$\bmatrix{#1}{#3}$\kern#2}
  \,\vbox{\offinterlineskip\hbox to\wd1{\hfil\copy1\hfil}
    \kern 3pt\hbox to\wd1{\rightarrowfill}}\,}

\def\leftarrowmat#1#2#3{
  \setbox1=\hbox{\small\kern#2$\bmatrix{#1}{#3}$\kern#2}
  \,\vbox{\offinterlineskip\hbox to\wd1{\hfil\copy1\hfil}
    \kern 3pt\hbox to\wd1{\leftarrowfill}}\,}

\def\rightarrowbox#1#2{
  \setbox1=\hbox{\kern#1\hbox{\small #2}\kern#1}
  \,\vbox{\offinterlineskip\hbox to\wd1{\hfil\copy1\hfil}
    \kern 3pt\hbox to\wd1{\rightarrowfill}}\,}

\def\leftarrowbox#1#2{
  \setbox1=\hbox{\kern#1\hbox{\small #2}\kern#1}
  \,\vbox{\offinterlineskip\hbox to\wd1{\hfil\copy1\hfil}
    \kern 3pt\hbox to\wd1{\leftarrowfill}}\,}


\newcount\countA
\newcount\countB
\newcount\countC

\def\monthname{\begingroup
  \ifcase\number\month
    \or January\or February\or March\or April\or May\or June\or
    July\or August\or September\or October\or November\or December\fi
\endgroup}

\def\dayname{\begingroup
  \countA=\number\day
  \countB=\number\year
  \advance\countA by 0 
  \advance\countA by \ifcase\month\or
    0\or 31\or 59\or 90\or 120\or 151\or
    181\or 212\or 243\or 273\or 304\or 334\fi
  \advance\countB by -1995
  \multiply\countB by 365
  \advance\countA by \countB
  \countB=\countA
  \divide\countB by 7
  \multiply\countB by 7
  \advance\countA by -\countB
  \advance\countA by 1
  \ifcase\countA\or Sunday\or Monday\or Tuesday\or Wednesday\or
    Thursday\or Friday\or Saturday\fi
\endgroup}

\def\timename{\begingroup
   \countA = \time
   \divide\countA by 60
   \countB = \countA
   \countC = \time
   \multiply\countA by 60
   \advance\countC by -\countA
   \ifnum\countC<10\toks1={0}\else\toks1={}\fi
   \ifnum\countB<12 \toks0={\sevenrm AM}
     \else\toks0={\sevenrm PM}\advance\countB by -12\fi
   \relax\ifnum\countB=0\countB=12\fi
   \hbox{\the\countB:\the\toks1 \the\countC \thinspace \the\toks0}
\endgroup}

\def\timestamp{\dayname, \the\day\ \monthname\ \the\year, \timename}

%

%


\def\frac#1#2{{#1 \over #2}}

\let\text\hbox
\def\overset#1\to#2{\ \mathop{\buildrel #1 \over #2}\ }
\def\underset#1\to#2{{\mathop{\buildrel #1 \over #2}}}

\def\myoperator{NOOP}
\def\mymathopsp{\futurelet\testchar\maybesubscriptop}
\def\maybesubscriptop{\ifx_\testchar \let\next\subscriptop
	\else \let\next\nosubscriptop \fi \next}
\def\subscriptop_#1{{\rm \myoperator}_{#1}}
\def\nosubscriptop{\futurelet\testchar\maybeparenop}
\def\maybeparenop{{\mathop{\rm \myoperator}\nolimits%
	\ifx(\testchar \relax\else\,\fi}}


\def\Ass{\mathop{\rm Ass}\nolimits}

\def\depth{\mathop{\rm depth}\nolimits\,}


\def\depth{\mathop{\rm depth}\nolimits}

\def\semidirect{{\mathop{\hbox to 1.1em{\hskip.5em\psline[linewidth=0.3pt](-.18,+.01)(.18,.15)(.18,.01)(-.18,.15)}}}}


\def\spm{\raise0.5pt\hbox{\bgroup\textfont2=\fivesy\relax$\pm$\egroup}}
\def\dotdiv{\mathop{\hbox{\rlap{\raise-0.1em\hbox{--}}\hskip0.1em\raise0.7ex\hbox{.}}\ }}


\def\uspace{\hbox{$\underline{\hbox to 1ex{\ \hfil}}$}}
\def\lineunderblank{\underline{\ \hbox to 0.5ex{$\vphantom{a}$\hfill\ }}}
\def\lineunderblankb#1{\underline{\ \hbox to #1{$\vphantom{a}$\hfill\ }}}
\def\lineaboveblank{\overline{\ \hbox to 0.5ex{$\vphantom{a}$\hfill\ }}}
\def\lineaboveblank{\overline{\ \hbox to 0.5ex{$\vphantom{a}$\hfill\ }}}
\def\blank#1{\hbox{$\underline{\hbox to #1em{\ \hfill}}$}}


\def\qedbox{\hbox{\vbox{\hrule\hbox{\vrule\kern3pt\vbox{\kern6pt}\kern3pt\vrule}\hrule}}}
\def\qed{\unskip\hfill\qedbox\vskip\thmskip} 
\def\eqed{\eqno\hbox{\quad\qedbox}} 


\def\today{\ifcase\month\or January\or February\or March\or
  April\or May\or June\or July\or August\or September\or
  October\or November\or December\fi
  \space\number\day, \number\year}\relax


\newwrite\isauxout
\openin1\jobname.aux
\ifeof1\message{No file \jobname.aux}
       \else\closein1\relax\input\jobname.aux
       \fi
\immediate\openout\isauxout=\jobname.aux
\immediate\write\isauxout{\relax}

\newwrite\iscontout
\immediate\openout\iscontout=\jobname.cont
\immediate\write\iscontout{\relax}

\def\contlist#1#2#3{\ifchap\vskip1pt\noindent \ \ \ifx0#1#2\hfill#3%
		\else\ref{#1}: #2\hfill#3%
                \fi%
	\else \vskip1pt\noindent\ref{#1}: #2\hfill#3\fi}
\def\contlist#1#2#3{\vskip1pt\noindent \ \hskip1em%
	\ifchap\ifx0#1#2\hfill#3%
		\else\ref{#1}: #2\hfill#3%
                \fi%
	\else \ref{#1}: #2\hfill#3\fi}
\def\contnosectlist#1#2#3{\vskip1pt\noindent\hskip1em#2\hfill#3}


%
\newwrite\isindexout
\immediate\openout\isindexout=\jobname.index
\font\footttt=cmtt6

\def\padno{
\ifnum\pageno<10 00\the\pageno\else\ifnum\pageno<100 0\the\pageno\else\the\pageno\fi\fi}

\def\index#1{\rlap{\ifdraft\draftind{{\footttt\Blue{#1\ \ }}}\unskip\fi%
\write\isindexout{\noexpand\indexentry {#1}{\padno}}}\unskip\ignorespaces}

\def\indexn#1{\rlap{\ifdraft\draftind{{\footttt\Brown{#1\ \ }}}\unskip\fi%
\write\isindexout{\noexpand\indexentryn {#1}{\padno}}}\unskip\ignorespaces}

\def\obeyspaces{\catcode`\ =\active\catcode`\	=\active}
{\obeyspaces\global\let =\space} 
{\obeyspaces\gdef {\hskip.5em}\gdef	{\hskip3em}}
{\catcode`\^^M=\active %
\gdef\obeylines{\catcode`\^^M=\active \gdef^^M{\vskip.0ex\ }}}%

\font\stt=cmtt10
\def\dMaccode{\vskip3ex\bgroup\advance\leftskip1ex\parindent=0pt%
	\baselineskip=9pt\obeylines\obeyspaces%
	\catcode`\^=11\catcode`\_=11\catcode`\#=11\catcode`\%=11%
	\catcode`\{=11\catcode`\}=11%
	\stt}
\def\enddMaccode{\removelastskip\unskip\egroup}
\def\Maccode{\bgroup \catcode`\^=11 \catcode`\_=11 \stt"}
\def\endMaccode{\unskip"\egroup}
\def\Maccode{\bgroup \catcode`\^=11 \catcode`\_=11 \stt}
\def\endMaccode{\unskip\egroup}

\gdef\poetry{\vskip0pt\bgroup\parindent=0pt\obeylines\obeyspaces}
\gdef\poetrycenter{\vskip0pt\bgroup\parskip=0pt\leftskip=0pt plus 1fil\rightskip=0pt plus
1fil\parfillskip=0pt\obeylines}
\gdef\endpoetry{\vskip0pt\egroup}

\magnification=1200
\abovedisplayskip=5pt
\belowdisplayskip=5pt

\newwrite\isauxout
\openin1\jobname.aux
\ifeof1\message{No file \jobname.aux}
       \else\closein1\relax\input\jobname.aux
       \fi
\immediate\openout\isauxout=\jobname.aux
\immediate\write\isauxout{\relax}

\def\uq{\underline q}

\def\Iuqm{I_{\underline q,m}}
\def\Ioneqm{I_{\{1\},m}}

\drafttrue
\draftfalse
\colortrue
\ifcolor\input colordvi\else\input blackdvi\fi%
\ifdraft\footline={{\rm \hfill page \folio -- \today\ --
\timehhmm\hfill}}\else
\ifnum\folio>1\footline={{\rm\hfill\folio\hfill}}\else\footline={{\rm\hfill\folio\hfill}}\fi\fi

\centerline{\sectionfont Predicted decay ideals}
\centerline{\sectionfont Sarah Jo Weinstein and Irena Swanson}
\centerline{Reed College, 3203 SE Woodstock Boulevard, Portland, OR 97202}
\centerline{\tt weinstein.sarah.j@gmail.com, iswanson@reed.edu}

\vskip 0.5cm
\bgroup
\narrower\narrower
\noindent
{\bf Abstract.}
\sl
For arbitrary positive integers $q_1 \ge q_2 \ge \cdots$
we construct a family of monomial ideals
such that for each positive integer $e$ and for each ideal $I$ in the family,
the number of associated primes of $I^e$ is $q_e$.
We present the associated primes explicitly.

\vskip 0.3cm
{\smallfont
\noindent
{\it Keywords and phrases:}
primary decomposition,
powers of ideals,
associated primes,
monomial ideals.

\noindent
{\it 2010 Mathematics Subject Classification. Primary 13A15, 13F20.}

}

\egroup 
\vskip 0.5cm

Given a sequence $q_1, q_2, \ldots$ of positive integers,
is it possible to find an ideal $I$ in a Noetherian ring
such that for each positive integer $e$,
the number of associated primes of~$I^e$ equals $q_e$?
If yes,
we say that such $I$ {\bf corresponds} to the sequence $\{q_e\}_{e \ge 1}$,
and that $\{q_e\}_{e \ge 1}$ is {\bf realizable} by~$I$.

The main result of this paper is
that every non-increasing sequence
is realizable.
Moreover,
for any non-increasing sequence we construct
an infinite family of corresponding monomial ideals,
all in the same polynomial ring over a field.
The construction
was the work of the first author's summer project
and subsequent senior thesis
under the supervision of the second author.
In general it is not an easy problem to describe
the associated primes of all powers of an ideal,
even of a monomial ideal,
but we do so for our constructed ideals.

By a result of Brodmann~\cite{Brodmann},
a realizable sequence is eventually constant.
We do not know if the eventually-constant restriction is sufficient for
realizability in general.
There exist monomial ideals $I$
for which the sequence $\{|\Ass(R/I^e)|\}_{e \ge 1}$ is not non-increasing,
which means that
our main result does not cover all monomial ideals.
For example,
the monomial ideal $(xy,xz,yz)$ has three associated primes
and its second power has four associated primes.
Similarly,
not every eventually stable sequence of positive integers
is realizable by a monomial ideal.
For example,
any monomial ideal corresponding to the sequence $\uq$ with $q_1 = 1$
must have all other $q_e$ equal to~$1$
because $q_1 = 1$ says that the ideal is primary,
and so by \ref{factsmonprim}~(2)
all the powers of the ideal are also primary.
Without the restriction to monomial ideals,
$q_1 = 1$ does allow arbitrarily large $q_2$;
the paper~\cite{KS} constructs prime ideals
in polynomial rings over fields
for which the number of associated primes of the second power
is exponential as a function in the number of variables in the ring.

Theorem 1 in the paper~\cite{BHH} by Bandari, Herzog and Hibi says
that for an arbitrary non-negative integer $n$
there exists a monomial ideal $I$ in a polynomial ring~$R$ in $4 + 2n$ variables
such that in the sequence $\{\depth(R/I^e)\}_{e \ge 1}$,
the first $n+1$ odd-numbered entries are~$0$,
the first $n$ even-numbered entries are~$1$,
and the rest of the entries are~$2$.
This means in particular that the maximal ideal
is associated to the first $n+1$ odd powers
and to no other powers of $I$.
Experimental work with Macaulay2~\cite{GS}
on these Bandari-Herzog-Hibi ideals
shows that for $n \ge 2$
the sequence $\{|\Ass(R/I^e)|\}_e$ is not monotone
and it does not have the same fluctuations
as the sequence $\{\depth(R/I^e)\}_{e \ge 1}$.

We thank Samuel Johnston, another Reed student,
for the conversations on the topic.

\section{Basics of primary decomposition and monomial ideals}

\facts[factsprimdec]
Let $R$ be a commutative ring with identity.
\ritem
For any ideal $I$ in $R$, any $x \in R$ and any non-negative integer $n$,
if $I:x^n = I:x^{n+1}$, then $I:x^n = I : x^m$ for all integers $m \ge n$.
We denote this stable colon ideal as $I:x^{\infty}$.
\ritem
For any ideal $I$ in $R$ and any $x \in R$,
if $I:x^n = I:x^{\infty}$, then $I = (I:x^n) \cap (I+(x^n))$.
\ritem
The previous part implies
that the set of associated primes of $I$
is either associated to $I : x^n$ or to $I + (x^n)$.
By properties of primary decompositions,
a prime ideal associated to $I : x^n$ is also associated to $I$,
but a prime ideal associated to $I + (x^n)$ need not be associated to $I$.
\ritem
If $P$ is a prime ideal minimal over $I$,
then $P$ is associated to all powers of~$I$.
\ritem
Let $R$ be a polynomial ring in variables $x_1, \ldots, x_n$ and $y_1, \ldots,
y_m$ over a field~$F$,
$I$ an ideal in $F[x_1, \ldots, x_n]$
and $J$ an ideal in $F[y_1, \ldots, y_m]$.
\riitem{(i)}
If $I = q_1 \cap \cdots \cap q_r$ is an irredundant primary decomposition,
then
$IR = (q_1 R) \cap \cdots \cap (q_r R)$ is an irredundant primary decomposition.
\riitem{(ii)}
$(IR)(JR) = (IR) \cap (JR)$.
\riitem{(iii)}
The intersection of irreducible primary decompositions of $IR$ and of $JR$
is an irreducible primary decompositions of $(IR)(JR)$.

\facts[factsmonprim]
Let $R$ be a polynomial ring in a finite number of variables over a field.
A~monomial ideal is an ideal generated by monomials,
and it always has a unique minimal monomial generating set.
The intersection of two monomial ideals
is the ideal generated by the least common multiples
of pairs of monomials where the two monomials are taken
from generating sets of the two ideals.

Let $I$ be a monomial ideal.
\ritem
$I$ is prime if and only if $I$ is generated by a subset of the variables.
\ritem 
$I$ is primary if and only if every variable
that divides a minimal monomial generator of $I$
is in $\sqrt{I}$.
\ritem
A prime ideal is associated to $I$ if and only if it is
a monomial ideal of the form $I:f$ for some monomial $f$
(that is necessarily not in $I$).
\ritem
Let $n$ be a non-negative integer.
If $I$ is a monomial ideal and if each variable dividing a monomial $m$
appears to exponent at most $n$ in each minimal generator of $I$,
then $I : m^n = I : m^\infty$
is the ideal obtained from $I$
by replacing (in each minimal generator of $I$)
each variable appearing in $m$ by $1$.
\ritem
A variable $x$ appears in a minimal generator of $I$
if and only if some associated prime ideal of $I$ contains $x$.

\lemma[lemmamultelt]
Let $I$ be a monomial ideal and $x,y$ variables.
Suppose that for any minimal generator $w$ in $I$,
if $x$ divides $w$ then ${yw \over x} \in I$.
Then for any prime ideal $P$ associated to a power of $I$,
if $x \in P$ then $y \in P$.
\endb

\proof
Let $e$ be a positive integer,
let $P$ be associated to $I^e$ and $x \in P$.
Write $P = I^e : f$ for some monomial~$f$ not in $I^e$.
Then
$xf = N_0 N_1 \cdots N_e$
for some minimal generators $N_1, \ldots, N_e$ of $I$
and some monomial $N_0$.
If $x$ divides $N_0$,
then $f \in (N_1 \cdots N_e) \in I^e$,
which is a contradiction.
Thus $x$ must divide $N_i$ for some $i > 0$.
Without loss of generality $x$ divides $N_1$.
By assumption then ${y N_1 \over x} \in I$,
so that $yf = N_0 {y N_1 \over x} N_2 \cdots N_s \in I^e$.
Thus $y \in I^e : f = P$.
\qed

\section{Construction and notation}[sectconstr]

Let $\uq = \{q_1,q_2,\ldots, q_{n-1}, q_n, q_n, q_n, \ldots\}$
be a non-increasing sequence of positive integers,
where $q_n = q_{n+1} = q_{n+2} = \cdots$.
For induction we allow the exceptional empty sequence $\uq = \{\}$
with $n = 0$.

We fix an arbitrary integer $m \ge n$.

Let $F$ be a field.
For $\uq = \{\}$
we set $\Iuqm$ to be the zero ideal in the ring $R_{\uq} = F$.
For non-empty $\uq$ we proceed as below.
We first measure the differences between adjacent entries in $\uq$:
$$
t_i(\uq) = \cases{
       q_{n}-1, & if $i = n$; \cr
       q_{i}-q_{i+1}-1, & if $i = 1, \ldots, n-1$. \cr
     }
$$
With these we next define sets $J_{\uq}$ and $K_{\uq}$,
the ring $R_{\uq}$,
the ideal $\Ioneqm$ and elements $M_j$ in~$R_{\uq}$.
The listed generators of the ring $R_{\uq}$ are variables over $F$.
$$
\eqalignno{
J_{\uq} &= \{ i \in \{ 1,\ldots,n-1 \} : t_i(\uq) \ge 0 \}, \cr
K_{\uq} &= \{ i \in \{ 1,\ldots,n \} : t_i(\uq) \ge 1 \}, \cr
R_{\uq} &= F[a,b,x_j,y_{kl} : j \in J_{\uq},k \in K_{\uq}, 1\le l \le
t_k(\uq)], \cr
\Ioneqm &= (a^{m+2}, a^{m+1}b,ab^{m+1}, b^{m+2}), \cr
M_j &= a^mb^{m-j+1}x_j \hbox{ for any $j \in J_{\uq}$}.  \cr
}
$$
By standard convention empty products equal 1.
For all $j \in J_{\uq} \cup \{n\}$,
we set
$$
Z_j = \prod_{l=1}^{t_j(\uq)}y_{jl},
\hskip 1em
Y_j = Z_n \cdot \!\!\!\!\!\! \prod_{s \in J_{\uq}, s \ge j}\!\!\!\!\!\! Z_s.
$$
Thus $Y_n = Z_n$,
$Z_j = 1$ for $j \in J_{\uq} \setminus K_{\uq}$,
$Z_j$ and $Z_{k}$ have greatest common divisor equal to~$1$
for distinct $j$ and $k$,
and $Y_j$ is a multiple of $Y_{k}$ if $j<k$.
Finally, we define the ideals
$$
\Iuqm
= \left(M_jY_{j} : j \in J_{\uq}\right)+\Ioneqm Y_n.
$$
We prove in this paper that the $\Iuqm$ realize~$\uq$.
By \ref{factsprimdec}~(5)(i),
without loss of generality we may add extra variables to $R_{\uq}$;
this is needed for induction.

When $\uq$ is the constant sequence $\{c\}$,
any positive integer $n$
can be used as the position beyond which the sequence stabilizes.
With a choice of $n$
the construction gives $J_{\uq} = \emptyset$,
and $t_n (\uq)= c-1$,
so that $K_{\uq} = \{n\}$ if $c > 1$ and $K_{\uq} = \emptyset$ if $c=1$.
Thus
$R_{\uq} = F[a,b, y_{n1}, y_{n2}, \ldots, y_{n,c-1}]$,
$\Iuqm
= (a^{m+2}, a^{m+1}b, a b^{m+1}, b^{m+2}) (y_{n1} y_{n2} \cdots y_{n,c-1})$.
By \refs{factsprimdec}~(5) and \refn{factsmonprim}~(2) and~(5),
for all positive integers $e$,
the set of associated primes of $\Iuqm^e$
equals the set $\{(a,b)$, $(y_{n1})$, $\ldots$, $(y_{n,c-1})\}$
of $c$ elements.
Thus $\Iuqm$ indeed realizes the constant sequence $\{c\}$.
Note that when $c=1$,
then $R_{\uq}=F[a,b]$ and $\Iuqm=\Ioneqm$,
so that the notation for $\Ioneqm$ makes sense.

\lemma[lmreduction]
a) $(a^{m+2} Y_n, a b^{m+1} Y_n) \cdot a^m b^2 Y_n \subseteq \Ioneqm^2$.

b)
Let $j, k \in J_{\uq}$ with $j > k$.
Then $x_j M_k Y_k \in (x_k Z_k M_j Y_j)$.
\endb

\proof
Part a) holds because
$a^{m+2} Y_n a^mb^2  Y_n =  (a^{m+1} bY_n)^2$
and
$a b^{m+1} Y_n a^mb^2 Y_n = (a^{m+1} b Y_n) (b^{m+2} Y_n)$.
Part~b) follows trivially.
\qed

\section{Associated primes of powers of $\Iuqm$}

Throughout we use the notation from \ref{sectconstr}.

\prop[propminpr]
Let $n$ and $e$ be positive integers.
Then the following statements hold.

a) $\{(a,b), (y_{nl}): 1 \le l \le t_n(\uq)\} \subseteq \Ass(R/\Iuqm^e)$.

b)
If $P \in \Ass(R/\Iuqm^e)$,
then either
$P = (y_{nl})$ for some $l \in \{1, \ldots, t_n(\uq)\}$
or else $(a,b)\subseteq P$
and no $y_{nl}$ is in $P$.
%
\endb

\proof
By definition
$\Iuqm = I_0 Z_n$ for some monomial ideal $I_0$
whose minimal generating set uses no variables $y_{nl}$.
Thus $\Iuqm^e = I_0^e Z_n^e$.
By \ref{factsprimdec}~(5),
$\Ass(R/\Iuqm^e) = \Ass(R/(Z_n)^e) \cup \Ass(R/I_0^e)$.
The prime ideal $(a,b)$ is the unique prime ideal minimal over $I_0^e$,
and the prime ideal $(y_{nl})$ is minimal over $(Z_n)^e$.
This in combination with \ref{factsprimdec}~(4) proves a).
The only associated primes of the principal ideal $(Z_n)^e$ are the $(y_{nl})$,
and all associated primes of $I_0^e$ contain $(a,b)$.
The latter prime ideals do not contain any~$y_{nl}$
by \ref{factsmonprim}~(5).
This proves b).
\qed

\prop[propxchain]
Let $P$ in $\Ass(R/\Iuqm^e)$
and let $x_k$ or $y_{kl}$ be in $P$.
If $j \in J_{\uq}$ and $j>k$, then $x_j \in P$.
\endb

\proof
Let $N$ be any monomial minimal generator of $\Iuqm$
that is a multiple of $w = x_k$ or of $w = y_{kl}$.
Then ${x_j N \over w} \in \Iuqm$ by the definition of the generators.
Thus by \ref{lemmamultelt},
$x_j \in P$.
\qed

\prop[propesmall]
Let $n$ and $e$ be positive integers,
let $P$ be in $\Ass(R/\Iuqm^e)$,
and let $x_k$ or $y_{kl}$ be in $P$ for some $k \in J_{\uq}$.
Then $k \ge e$.
\endb

\proof
Suppose for contradiction that $k < e$.
Then $e \ge k+1 \ge 2$.
We know that $P = \Iuqm^e : f$ for some monomial $f$
necessarily not in $\Iuqm^e$.
By assumption $w = x_k \in P$ or $w = y_{kl} \in P$.
Thus $w f = N_0 N_1 \cdots N_e$
for some minimal generators $N_1, \ldots, N_e$ of $\Iuqm$
and some monomial $N_0$.

We next prove that $f \in \Iuqm^e$,
which gives a contradiction.
If the variable $w$ divides $N_0$,
then $f \in (N_1 \cdots N_e) \in \Iuqm^e$.
Thus we may assume that $w$ divides $N_i$ for some $i > 0$.
Without loss of generality $w$ divides~$N_1$,
so that $N_1 = a^m b^{m-i+1} x_i Y_i$ for some $i \le k$
and $N_1 \over w$ is a multiple of $a^m b^{m-k+1} Y_n$.
It remains to prove that $g = a^m b^{m-k+1}Y_n N_2 \cdots N_e$
is an element of~$\Iuqm^e$.
Recall that $e \ge 2$
and that $j \le n-1 \le m-1$ for $j \in J_{\uq}$,
so that $m-j+1 \ge~2$.
By \ref{lmreduction}~a) it suffices to prove that
$I' = a^m b^{m-k+1} Y_n
(a^{m+1}b Y_n, b^{m+2} Y_n, a^m b^2 Y_j: j \in J_{\uq})^{e-1}$
is contained in $\Iuqm^e$.
An arbitrary monomial in $I'$ can be written as
a multiple of
$a^mb^{m-k+1} Y_n(a^{m+1}b Y_n)^{c_1}(b^{m+2} Y_n)^{c_2 }(a^mb^2Y_n)^{c_3}$,
where $c_1, c_2, c_3$ are non-negative integers and
$c_1 + c_2  + c_3 + 1 = e$.
This monomial can be rewritten as
$a^{m-e} b^{e-k-1}
\cdot (a^{m+2}Y_n)^{c_1}
(ab^{m+1}Y_n)^{c_2  + 1}
(a^{m+1}bY_n)^{c_3}$,
which is an element of $\Iuqm^{c_1 + c_2 +1 + c_3} = \Iuqm^e$.
\qed

\prop[propbige]
For all $e \ge n \ge 1$,
$\displaystyle
\Ass(R/\Iuqm^e) = \{(a,b),(y_{nl}) : 1 \le l \le t_n(\uq)\}$.
\endb

\proof
Let $P \in \Ass(R/\Iuqm^n)$.
By \ref{propesmall},
$P$ does not contain any $x_k$ or $y_{kl}$ for $k \in J_{\uq}$.
Thus $P$ can only contain variables $a, b$ or the $y_{nl}$.
By \ref{propminpr}~b),
$P$ can only be one of $(y_{nl})$ or $(a,b)$.
By \ref{propminpr}~a),
all these ideals are associated.
\qed

\prop[propxy]
Let $e$ be a positive integer and let $P$ be associated to $\Iuqm^e$.
Let $k,j \in J_{\uq}$ with $j \ge k$.
Then $P$ contains at most one of $x_k, y_{kl}, y_{jp}$.
\endb


\proof
Let $x$ stand for $x_k$ or for $y_{kl}$,
and let $y$ stand for $y_{jp} \not = x$,
or alternatively let $x = x_k$ and $y = y_{kl}$.
Suppose for contradiction that $x, y \in P$.
By \ref{propesmall} we know that $e \le k$,
and so $e < n \le m$.
Let
$$
\eqalignno{
I_1 &= \Ioneqm Y_n + (M_i Y_i : i \in J_{\uq}, i > j), \cr
I_2 &= \cases{
{1 \over y} (M_i Y_i : i \in J_{\uq}, i \not = k, i \le j), & if $x = x_k$; \cr
{1 \over y} (M_i Y_i : i \in J_{\uq}, k < i \le j), & if $x = y_{kl}$,  \cr
} \cr
I_3 &= \cases{
{1 \over xy} (M_k Y_k), & if $x = x_k$; \cr
{1 \over xy} (M_i Y_i : i \in J_{\uq}, i \le k), & if $x = y_{kl}$. \cr
} \cr
}
$$
Then $I_1, I_2, I_3$ are ideals in $R_{\uq}$
and $\Iuqm = I_1 + I_2 y + I_3 xy$.
We set $I_x = I_1 + I_2 + I_3 x$ and $I_y = I_1 + I_2 y + I_3 y$.

Suppose that $\Iuqm^e = I_x^e \cap I_y^e$.
Then the intersection of irredundant primary components of $I_x^e$ and $I_y^e$
is a possibly redundant primary decomposition of $\Iuqm^e$.
Thus any prime ideal associated to $\Iuqm^e$
is either associated to $I_x^e$ or to $I_y^e$.
Since the generators of $I_x$ involve no~$y$
and the generators of $I_y$ involve no $x$,
by \ref{factsmonprim}~(5),
none of these prime ideals contains both $x$ and $y$.

Thus it suffices to prove that $I_x^e \cap I_y^e = \Iuqm^e$.
Since $\Iuqm \subseteq I_x \cap I_y$,
it suffices to prove that $I_x^e \cap I_y^e \subseteq \Iuqm^e$.
We proceed by induction on $e$.

Common-factor observation:
For any $f \in I_1$,
$(f) I_x^{e-1} \cap (f) I_y^{e-1} \subseteq (f) (I_x^{e-1} \cap I_y^{e-1})$
is contained in $\Iuqm^e$ by induction on~$e$.

Since intersection of monomial ideals distributes over sums,
it suffices to prove that for any non-negative integers $u_1, u_2, u_3, v_1,
v_2, v_3$,
if $u_1 + u_2 + u_3 = e = v_1 + v_2 + v_3$,
then
$I_1^{u_1} I_2^{u_2} (I_3 x)^{u_3} \cap
I_1^{v_1} (I_2y)^{v_2} (I_3 y)^{v_3} \subseteq \Iuqm^e$.
Since $I_1, I_2, I_3$ use no variables $x, y$,
the intersection equals
$I_1^{u_1} I_2^{u_2} (I_3 x)^{u_3} y^{v_2 + v_3} \cap
I_1^{v_1} (I_2y)^{v_2} (I_3 y)^{v_3} x^{u_3}$.
If $u_3 \ge v_3$,
then the intersection is contained in 
$I_1^{v_1} (I_2y)^{v_2} (I_3 xy)^{v_3} \subseteq \Iuqm^e$,
and if $v_2 + v_3 \ge u_2 + u_3$,
then the intersection is contained in 
$I_1^{u_1} (I_2y)^{u_2} (I_3 xy)^{u_3} \subseteq \Iuqm^e$.
Thus we may assume that 
$u_3 < v_3$ and $v_2 + v_3 < u_2 + u_3$.
This in particular implies that it suffices to assume that
$v_3$ and $u_2$ are positive,
and even that $u_2 \ge v_2 + 2$.
This proves the base case $e = 1$ trivially.
So let $e > 1$.

The strategy for the rest of the proof is to keep lowering $v_3$ or $u_2$
at the expense of increasing $v_1$ or $u_1$, respectively,
thus using induction on $u_2 + v_3$ within induction on $e$.


It remains to prove that any monomial $w$ in
$I_1^{u_1} I_2^{u_2} (I_3 x)^{u_3} \cap I_1^{v_1} (I_2y)^{v_2} (I_3 y)^{v_3}$
is in~$\Iuqm^e$.
By assumption we can write $w$ as a product
of $u_1$ factors from $I_1$,
$u_2$ factors from $I_2$ and $u_3$ factors from $I_3 x$,
and simultaneously we can write $w$ as a product
of $v_1$ factors from $I_1$,
$v_2$ factors from $I_2y$ and $v_3$ factors from $I_3 y$.

If one of the two factorizations of $w$ as described in the previous paragraph
has $M_i Y_i$ as a factor from $I_1$,
then since the polynomial ring $R_{\uq}$ is a unique factorization domain,
the other factorization must have $x_i$ as a factor.
Since $x_i$ only appears as a factor in the minimal generator $M_i Y_i$,
and since by the common-factor observation
the other factorization does not have $M_i Y_i$ as one of the $e$ factors,
necessarily $x_i$ is a redundant factor in the other factorization.
Thus $w \in x_i I_1^{u_1} I_2^{u_2} (I_3 x)^{u_3}$
or $w \in x_i I_1^{v_1} (I_2y)^{v_2} (I_3 y)^{v_3}$.
Using that $u_2 v_3 > 0$ and \ref{lmreduction}~b),
we can rewrite $w$ with lower $u_2 + v_3$.
Thus we may assume that
$w \in \Ioneqm^{u_1} I_2^{u_2} (I_3 x)^{u_3} \cap
\Ioneqm^{v_1} (I_2y)^{v_2} (I_3 y)^{v_3}$.
By applying \ref{lmreduction}~a)
we may also reduce $u_2 + v_3$
in case $a^{m+2} Y_n$ or $a b^{m+1}Y_n$ are used as factors in $\Ioneqm$,
and thus we may assume that
$w \in (a^{m+1}bY_n, b^{m+2}Y_n)^{u_1} I_2^{u_2} (I_3 x)^{u_3} \cap
(a^{m+1}bY_n, b^{m+2}Y_n)^{v_1} (I_2y)^{v_2} (I_3 y)^{v_3}$.
By the common-factor observation,
$a^{m+1} b Y_n$ and $b^{m+2} Y_n$ are each a factor
in at most one factorization of $w$.
If $a^{m+1} b Y_n$ is a factor in one,
call $L$ the other intersectand ideal.
Since $e < m$,
by the $a$-degree count,
$w \in a L$.
Since $u_2 v_3 > 0$ and since $a M_j Y_j \in (a^{m+1}b Y_n)$ for all $j \in
J_{\uq}$,
we may rewrite $w$ with lower $u_2 + v_3$,
and we are done by induction.
If instead $b^{m+2} Y_n$ is a factor,
then again $w \in a L$ by the $a$-degree count
and we are done.
Thus we may assume that $u_1 = v_1 = 0$,
and so $u_3 = e - u_2$, $v_3 = e - v_2$, and
$$
\eqalignno{
w &\in
I_2^{u_2} (I_3 x)^{e-u_2} \cap (I_2y)^{v_2} (I_3 y)^{e-v_2} \cr
&\subseteq I_2^{u_2} (I_3 x)^{e-u_2} \cap (y)^e\cr
&= I_2^{u_2} (I_3 x)^{e-u_2} (y)^e\cr
&= (I_2y)^{u_2} (I_3 xy)^{e-u_2}\cr
&\subseteq \Iuqm^e.
&\qedbox
}
$$

\prop[propbigass]
For any positive integers $n$ and $e$,
the set of associated primes of $\Iuqm^e$
is a subset of
$$
\eqalignno{
T &= \{(a, b), (y_{nl}) : 1 \le l \le t_n(\uq)\} \cr
&\cup \bigcup_{k \in J_{\uq},k \ge e}
\left(
\{(a,b,x_j : j \in J_{\uq}, k \le j)\}
\cup
\{(a,b,y_{kl},x_j : j \in J_{\uq}, k < j) : l = 1, \ldots, t_k(\uq)\}
\right).
\cr
}
$$
\endb

\proof
Let $P \in \Ass(R/\Iuqm^e)$.
Since $\Iuqm$ is a monomial ideal,
$P$ is generated by a subset of the variables.
By \ref{propminpr}~b) and \ref{propxy},
$P$ contains at most one~$y_{kl}$.
By \ref{propminpr},
each $y_{nl}$ appears only in the associated prime $(y_{nl})$,
and all other prime ideals contain $(a,b)$.
Combination with \ref{propxy} says that if $y_{jl} \in P$,
then $x_k \not \in P$ for all $k \le j$.
If $y_{kl}$ or $x_k$ is in~$P$,
then by \ref{propxchain},
$x_j$ is in~$P$ for all $j$ in $J_{\uq}$ with $j > k$,
and by \ref{propesmall},
$k \ge e$.
\qed

\prop[propbasecase]
Let $n$ and $e$ be positive integers
and let $\uq$ be such that $|K_{\uq}| = 0$.
Then $\Ass(R/\Iuqm^e)$ is the set $T$ in \ref{propbigass},
with $q_e$ elements.
\endb

\proof
By assumption $K_{\uq} = \emptyset$,
the tail-end repeating entry in $\uq$ is $1$
and the differences of consecutive terms in $\uq$ are at most $1$.
Thus we can write
$$
\uq = \{r,\ldots, r,r-1,\ldots,r-1,\ldots,2,\ldots,2,1,\ldots,1, \ldots\},
$$
where each $i = 1, 2, \ldots, r = q_1$ appears $l_i$ times with $l_i \ge 1$
and we may take the count $l_1$ of the trailing $1$s to be any positive integer.
Set $j_k = \sum_{i=r-k+1}^r l_i$.
So $j_1 = l_r$, $j_2 = l_r + l_{r-1}$, and so on.
Then
$$
\eqalignno{
j_1 &< j_2 < \cdots < j_r = n \le m, \cr
t_i(\uq) &= \cases{
       0, & if $i \in \{j_1, \ldots, j_r\}$; \cr
       -1, & otherwise. \cr
     } \cr
J_{\uq} &= \{j_1, j_2, \ldots, j_{r-1}\}, \cr
R_{\uq} &= F[a,b,x_{j_1}, \ldots, x_{j_{r-1}}], \cr
\Iuqm &=
(a^{m+2},a^{m+1}b,ab^{m+1}b^{m+2},
a^m b^{m-j_1+1} x_{j_1}, \ldots, a^m b^{m-j_{r-1}+1} x_{j_{r-1}}).
}
$$
Thus by \ref{propbigass},
$$
\Ass(R/\Iuqm^e) \subseteq
\{(a,b)\}
\cup \bigcup_{s = 1, j_s \ge e}^{r-1}
\{(a,b,x_{j_s}, x_{j_{s+1}}, \ldots, x_{j_{r-1}})\}.
$$
If we prove that the inclusion is an equality,
then the count follows as well.
\ref{propbige} proves the case $e \ge n$,
so we may assume that $e < n$.
By \ref{propminpr}~a)
it remains to prove that
$P_s = (a, b, x_{j_s}, x_{j_{s+1}}, \ldots, x_{j_{r-1}})$
is associated to $\Iuqm^e$
for $e = 1, 2, \ldots, j_s$.
By \ref{factsmonprim}~(3)
it suffices to prove that $P_s = \Iuqm^e :~w_e$,
where
$$
w_e = \cases{
a^m b^{m + (e-1)(m+2)}, & if $s = 1$ and $e = 1, \ldots, j_{s-1} + 1$; \cr
a^m b^{m-j_{s-1} + (e-1)(m+2)} x_{j_{s-1}}, & if $s > 1$ and $e = 1, \ldots, j_{s-1} + 1$; \cr
a^m b^{e(m+1) - 1}, & if $e = j_{s-1} + 2, \ldots, j_s$, \cr
}
$$
where we use the convention that $j_0 = 0$ and $x_0 = 1$.

First let $e = 1$.
$$
\Iuqm : a^m b^{m-j_s+1} =
(a^2,a, b^{j_s}, b^{j_s+1},
b^{j_s  - j_1} x_{j_1},
\ldots,
b^{j_s  - j_{s-1}} x_{j_{s-1}},
x_{j_s},
x_{j_{s+1}}, \ldots, x_{j_{r-1}}).
$$
In case
$s = 1$,
we colon this with $b^{j_1 -1}$,
and
if $s > 1$
we colon with $b^{j_s - j_{s-1}-1}x_{s-1}$,
and in both cases we get $P_s$ as the result.
Thus $P_s = \Iuqm : w_1$
for $w_1$ equal to $a^m b^{m- j_{s-1}} x_{j_{s-1}}$
and $a^m b^m$, respectively.
This proves the case $e = 1$.

Now let $e \in \{2, \ldots, j_{s-1}+1\}$.
Note that $w_e = w_{e-1} b^{m+2}$,
and by induction on $e$ we have that $\Iuqm^{e-1} : w_{e-1} = P_s$.
Since $b^{m+2} \in \Iuqm$,
it follows that
$P_s = \Iuqm^{e-1} : w_{e-1} \subseteq \Iuqm^{e} : w_{e-1} \Iuqm \subseteq
\Iuqm^{e} : w_{e-1} b^{m+2}$
$= \Iuqm^{e} : w_e$.
Suppose for contradiction that $P_s$ is strictly contained in $\Iuqm^e : w_e$.
Then
$\Iuqm^e : w_e$ contains a monomial $M$
using variables $x_{j_1}, \ldots, x_{j_{s-1}}$ only.
By construction,
$M w_e$ has $a$-degree $m$
and $x_{j_i}$-degree~$0$ for all $i \ge s$,
so that by looking at the generators of $\Iuqm$,
we conclude that
$M w_e$ is contained in the following subideal of $\Iuqm^e$:
$$
M w_e \in (ab^{m+1}, b^{m+2}, a^m b^{m-j_1+1} x_{j_1}, \ldots, a^m b^{m-j_{s-1}+1} x_{j_{s-1}})^e.
$$
Thus there exist non-negative integers $c_1, c_2, c_3$
and a monomial $N$
such that $c_1 + c_2 + c_3 = e$
and
$M w_e
= N (ab^{m+1})^{c_1} (b^{m+2})^{c_2} (a^m b^{m-j_{s-1}+1})^{c_3}$.
By looking at the exponents of $a$ and~$b$ in the equation we get that
$$
\eqalignno{
m &\ge c_1 + c_3 m, \cr
m - j_{s-1} + (e-1)(m+2) &\ge c_1 (m+1) + c_2 (m+2) + c_3 (m-j_{s-1} + 1).
\cr
}
$$
The first inequality implies that $c_3 \le 1$.
If $c_3 = 1$, then $c_1 = 0$ and $c_2 = e-c_1-c_3=e-1$,
and the second inequality says that
$m - j_{s-1} + (e-1)(m+2) \ge (e-1) (m+2) + (m-j_{s-1} + 1)$,
which is impossible.
Thus $c_3 = 0$.
Then $c_2 = e - c_1$,
and the second inequality says that
$$
m - j_{s-1} + (e-1)(m+2)
\ge c_1 (m+1) + (e- c_1) (m+2)
= e (m+2) - c_1,
$$
so that
$- j_{s-1} \ge 2 -c_1$.
Hence $2 + j_{s-1} \le c_1 \le e \le j_{s-1} + 1$,
which is a contradiction.
This proves that $P_s$ is associated to $\Iuqm^e$
for $e = 1, \ldots, j_{s-1} + 1$.

It remains to prove that $P_s = \Iuqm^e : w_e$
for $e = j_{s-1}+2, \ldots, j_s$.
Note that
$$
\eqalignno{
a w_e &= a a^m b^{e(m+1)-1} = b^{m-e} (a^{m+1}b) (b^{m+2})^{e-1} \in \Iuqm^e, \cr
b w_e &=
b a^m b^{e(m+1)-1}
= a^{m-e} (ab^{m+1})^e \in \Iuqm^e, \cr
x_{j_i} w_e &=
x_{j_i} a^m b^{e(m+1)-1}
= b^{j_i - e} (b^{m+2})^{e-1} (a^m b^{m-j_i+1} x_{j_i}) \in \Iuqm^e
\hbox{ for all $j_i \ge e$}.
 \cr
}
$$
This proves that $P_s \subseteq \Iuqm^e : w_e$.
Suppose for contradiction that $\Iuqm^e : w_e$ is strictly bigger than $P_s$,
so that it contains a monomial in $x_{j_1}, \ldots, x_{j_{s-1}}$.
By~\ref{propbigass},
these variables are not in any prime ideal associated to $\Iuqm^e$,
so that
$1 \in \Iuqm^e : w_e$,
i.e., that $w_e \in \Iuqm^e$.
By looking at the $x_i$-degrees
and the $a$-degree of the generators of $\Iuqm$,
we get that
$a^m b^{e(m+1)-1} =$
$w_e \in (a b^{m+1}, b^{m+2})^e \subseteq (b^{m+1})^e$,
which is not possible.
Thus $P_s = \Iuqm^e : w_e$,
which proves that $P_s$ is associated to $\Iuqm^e$.
\qed

\prop[propinduction]%
\thmparens{Inductive Lemma}
Let $n$ be a positive integer.
For 
$k \in J_{\uq} \cup \{n\}$
let
$$
\eqalignno{
h_{k}(\uq) &= \{q_1- t_k(\uq), q_2-t_k(\uq), \ldots, q_k - t_k(\uq),
q_{k+1},q_{k+2}, \ldots, q_n\}, \cr
g_n(\uq) &= \{\}, \cr
g_{k}(\uq) &= \{\underbrace{q_{k+1},\ldots,q_{k+1}}_{k+1 \hbox{ times } },
q_{k+2},\ldots,q_n\} \hbox{ if $k < n$}. \cr
}
$$
Then the newly defined sequences are non-increasing.
The rings for the defining ideals for the new sequences
as constructed in \ref{sectconstr}
are generated over $F$ by a subset of the variables that define $R_{\uq}$.
For any integer $m \ge n$ we have the following equalities:

a) $\Iuqm :Z_k = \Iuqm :Z_k^\infty = I_{h_k(\uq),m} R_{\uq}$.

b) $J_{h_k(\uq)} = J_{\uq}$
and $K_{h_k(\uq)} \subseteq K_{\uq}$.

c) $\Iuqm + (Z_k) = I_{g_k(\uq),m} R_{\uq} + (Z_k) R_{\uq}$.

d) If $k \in K_{\uq}$,
then $|K_{h_k(\uq)}|,
|K_{g_k(\uq)}| <
|K_{\uq}|$.
\endb

\proof
It is clear that all $h_k(\uq)$ and
all $g_k(\uq)$ are non-increasing.

We compute
$$
\eqalignno{
t_i(h_n(\uq)) &= \cases{
q_n - t_n(\uq) - 1 = t_n(\uq) - t_n(\uq) = 0, & if $i = n$; \cr
q_i - t_n(\uq) - (q_{i+1} - t_n(\uq)) -1 = q_i - q_{i+1} - 1 = t_i(\uq), & if $1 \le i < n$, \cr
} \cr
\noalign{\noindent and for $k < n$,}
t_i(h_k(\uq)) &= \cases{
q_i - t_k(\uq) - (q_{i+1} - t_k(\uq)) - 1 = q_i - q_{i+1} - 1 = t_i(\uq), & if $i <
k$; \cr
q_{k} - t_{k}(\uq) - q_{k+1} - 1 = t_{k}(\uq) - t_{k}(\uq) = 0, & if $i = k$; \cr
q_i - q_{i+1} - 1 = t_i(\uq), & if $k < i < n$; \cr
q_n - 1 = t_n(\uq), & if $i = n$. \cr
} \cr
}
$$
Thus for all $k \in J_{\uq} \cup \{n\}$,
$J_{h_k(\uq)} = J_{\uq}$ and
$K_{h_k(\uq)} = K_{\uq} \setminus \{k\}$.
This proves b) and the first part of d).
Furthermore,
this means that $I_{h_k(\uq),m}$ is defined just like $\Iuqm$
but with replacing $Z_k$ with $1$.
Since all $y_{kl}$ appear in the generating set with exponent at most~$1$,
part a) follows by \ref{factsmonprim}~(4).

Now we prove c).
Since $\Iuqm \subseteq (Z_n)$ and $I_{g_n(\uq),m} = (0)$,
the proposition follows for $k = n$.
So let $k < n$.
If $t_k(\uq) = 0$,
then the ideal $(Z_k)$ is the whole ring so that c) follows trivially.
So we may suppose that $t_k(\uq) > 0$.
Then
$$
I_{\uq,m} + (Z_k)
= I_{\{1\},m} Y_n + (M_j Y_j : j \in J_{\uq}, j > k) + (Z_k).
$$
From
$$
\eqalignno{
t_i(g_{k}(\uq)) &= \cases{
-1, & if $i \le k$; \cr
t_i(\uq), & if $i > k$, \cr
} \cr
J_{g_k(\uq)} &= J_{\uq} \setminus \{1,2,\ldots, k\}, \cr
K_{g_k(\uq)} &= K_{\uq} \setminus \{1,2,\ldots, k\}. \cr
}
$$
we get that
$I_{g_k(\uq),m} = (M_j Y_j : j \in J_{\uq}, j > k) + I_{\{1\},m} Y_n$
and so
$I_{g_k(\uq),m} R_{\uq} + (Z_k) R_{\uq} = I_{\uq,m} + (Z_k)$.
This proves c)
and finishes the proof of d).
\qed

\thm[thmassprimes]
We use notation from \ref{sectconstr}.
For any positive integers $m \ge n$ and $e$,
any for any non-increasing positive-integer sequence $\uq=
\{q_1,q_2,\ldots,q_n, q_n, \ldots\}$,
the set of associated primes of $\Iuqm^e$
equals
$$
\eqalignno{
\{(a&,b), (y_{nl}) : 1 \le l \le t_n(\uq)\} \cr
&\cup \bigcup_{j \in J_{\uq},j \ge e}
\left(
\{(a,b,x_i : i \in J_{\uq}, j \le i)\}
\cup
\{(a,b,y_{jl},x_i : i \in J_{\uq}, j < i) : l = 1, \ldots, t_j(\uq)\}
\right).
\cr
}
$$
\endb

\proof
By~\refs{propminpr} and \refn{propbigass}
it suffices to prove that the set $S$
in the second row is a subset of $\Ass(R/\Iuqm^e)$.

The proof is by induction on $|K_{\uq}|$.
The base case $|K_{\uq}|=0$ is proved in~\ref{propbasecase},
so we may assume that $K_{\uq}$ is not empty.
In particular,
variables $y_{kl}$ appear in the definition of $\Iuqm$
for some $k \in K_{\uq}$ and all $l = 1, \ldots, t_k$ with $t_k > 0$.
Assume that the theorem holds for any non-increasing sequence~$\underline p$
such that $|K_{\underline p}| < |K_{\uq}|$.
In particular,
with notation as in \ref{propinduction},
the theorem holds for the sequences $h_k(\uq)$ and $g_k(\uq)$.
Since every $y_{kl}$
appears to exponent at most~$1$ in a minimal generating set of $\Iuqm$,
we have that
$I_{h_k(\uq),m}^e = (\Iuqm : Z_k^\infty)^e
= \Iuqm^e : Z_k^\infty$.
Every associated prime ideal of $\Iuqm^e : Z_k^\infty$
is associated to $\Iuqm^e$.
Thus by induction the following ideals are associated to
$I_{h_k(\uq),m}^e$ and therefore to $\Iuqm^e$:
$$
\bigcup_{j \in J_{\uq},j \ge e} \!\!\!
\{(a,b,x_i : i \in J_{\uq}, j \le i)\}
\, \cup \!\!\!
\bigcup_{j \in J_{\uq},j \ge e, j \not = k} \!\!\!
\{(a,b,y_{jl},x_i : i \in J_{\uq}, j < i) : l = 1, \ldots, t_j(\uq)\}.
$$
If $k = n$,
this set is the same as $S$.
So we may assume that $k < n$.
It remains to prove that for each $l = 1, \ldots, t_k(\uq)$,
$(a,b,y_{kl},x_i : i \in J_{\uq}, k < i)$ is associated to $\Iuqm^e$.
Since we just established that
$(a,b,x_i : i \in J_{\uq}, k \le i)$ is associated to~$\Iuqm^e$,
it means by \ref{factsmonprim}~(5)
that $x_k$ appears in a minimal generator of~$\Iuqm^e$.
But the only generator of~$\Iuqm$ in which $x_k$ appears
is also a multiple of $y_{kl}$ (for $l = 1, \ldots, t_k(\uq)$).
This means that $y_{kl}$ appears in a minimal generator of $\Iuqm^e$.
Then again by \ref{factsmonprim}~(5),
$y_{kl}$ must appear in some associated prime of $\Iuqm^e$.
By \ref{propbigass},
we have only one option for such a prime ideal,
and it is the one in $S$.
This finishes the proof of the theorem.
\qed

\thm\label{thm:neat}
Let $\uq$ be a non-increasing sequence of positive integers
with $q_n = q_{n+1} = \ldots$.
Then there exists a polynomial ring $R_{\uq}$
finitely generated over a field by at most $2 + q_1$ variables
such that for every integer $m \ge n$
there exists a monomial ideal $\Iuqm$ in $R_{\uq}$
such that for all positive integers~$e$,
the number of prime ideals associated to $R_{\uq}/\Iuqm^e$ is $q_e$.
\endb

\proof
We set $R_{\uq}$ and $\Iuqm$ as in \ref{thmassprimes}.
First let $e \ge n$.
By \ref{propbige} (or by \ref{thmassprimes}),
$\Ass(R/\Iuqm^e) = \{(a,b),(y_{nl}) : 1 \le l \le t_n(\uq)\}$,
and the size of this set is $1 + t_n(\uq) = q_n = q_e$.
Now suppose that $e < n$.
Then the number of associated primes of $\Ass(R/\Iuqm^e)$
listed in \ref{thmassprimes} is:
$$
1 + t_n(\uq)
+ \sum_{j \in J_{\uq}, j \ge e} \left(1 + t_j(\uq) \right)
= q_n
+ \sum_{j \in J_{\uq}, j \ge e} (q_j - q_{j+1})
= q_e. \eqed
$$

\bigskip\bigskip
\leftline{\bf References}
\bigskip

\baselineskip=9.9pt
\parindent=3.6em
\setbox1=\hbox{[999]} 
\newdimen\labelwidth \labelwidth=\wd1 \advance\labelwidth by 2.5em
\ifnumberbibl\advance\labelwidth by -2em\fi

\thmno=0

\bitem{BHH}
S. Bandari, J. Herzog, and T. Hibi,
Monomial ideals whose depth function
has any given number of strict local maxima,
{\it Ark. Mat.} {\bf 52} (2014), 11--19.

\bitem{Brodmann}
M. Brodmann,
Asymptotic stability of Ass($M/I^nM$),
{\it Proc. Amer. Math. Soc.} {\bf 74} (1979), 16--18.

\bitem{GS}
D. Grayson and M. Stillman,
Macaulay2, a software system for research in algebraic geometry,
available at {\tt http://www.math.uiuc.edu/Macaulay2}.

\bitem{KS}
J. Kim and I. Swanson,
Many associated primes of powers of primes,
arXiv:1803.05456,
preprint.
Published online in {\it J. Pure and Appl. Algebra} in February 2019.

\end

\example
Let $\uq = \{6,5,5,4,2,1\}$. Then
$$
\Iuqm = (a^8,a^7b,ab^7,b^8,a^6b^2x_5,a^6b^3x_4y_{4,1},a^6b^4x_3y_{4,1},a^6b^6x_1y_{4,1})
$$
and the sets of associated primes for the powers are as below:
$$
\eqalignno{
\Iuqm &:
\{(a,b),
(a,b,x_5),
(a,b,x_5,x_4),
(a,b,x_5,y_{4,1}),
(a,b,x_5,x_4,x_3),
(a,b,x_5,x_4,x_3,x_1)\}, \cr
\Iuqm^2 &:
\{(a,b),(a,b,x_5),(a,b,x_5,x_4),(a,b,x_5,y_{4,1}),(a,b,x_5,x_4,x_3)\}, \cr
\Iuqm^3 &:
\{(a,b),(a,b,x_5),(a,b,x_5,x_4),(a,b,x_5,y_{4,1}),(a,b,x_5,x_4,x_3)\}, \cr
\Iuqm^4 &:
\{(a,b),(a,b,x_5),(a,b,x_5,x_4),(a,b,x_5,y_{4,1})\}, \cr
\Iuqm^5 &:
\{(a,b),(a,b,x_5)\}, \cr
\Iuqm^e &:
\{(a,b)\} \hbox{ for all $e \ge 6$}. \cr
}
$$

\end

\lemma[lemmamultcolon]
Let $I$ be a monomial ideal and $x,y$ variables.
Suppose that for any minimal generator $w$ in $I$
and for any positive integer $c$,
$x^c$ divides $w$ if and only if $y^c$ divides $w$.
Then $(x,y)$ is not contained in any associated prime ideal of $I$.
\endb

\proof
Note that $I + (xy) = I_0 + (xy)$
for some monomial ideal $I_0$ in which $x$ does not divide
any minimal generator.

Below we use the fact that
For any $x \in R$ and any ideal $I$ in $R$,
$\Ass(R/I) \subseteq \Ass(R/(I + (x)) \cup \Ass(R/(I : x))$.

Let $P$ be associated to $I$ and $(x,y) \subseteq P$.
If $P$ is associated to $I + (xy)$,
then by the fact,
$P$ is associated to $(I + (xy)) + (y) = I_0 + (y)$
or to $(I + (xy)) : y = I_0 + (x)$.
But the minimal generators of $I_0 + (y)$ do not involve~$x$
and the minimal generators of $I_0 + (x)$ do not involve~$y$,
so that by \ref{factsmonprim}~(5),
$P$ is not associated to $I + (xy)$.
Thus by the fact,
$P$ is associated to $I : (xy)$.
But $I : (xy)$ has the same defining property as $I$,
so that by repetition of the argument
$P$ is associated to $I : (xy)^c$ for all positive integers~$c$.
In particular,
$P$ is associated to $I : (xy)^\infty$.
But $x$ and $y$ do not appear in any minimal generator of $I : (xy)^\infty$,
so that by \ref{factsmonprim}~(5),
$P$ cannot be associated to it.
Thus we get a contradiction.
\qed

We prove next that $\Iuqm^e = I_x^e \cap I_y^e$.
If $e = 1$,
the intersection equals
$I_1 + I_2 y + I_3 xy = \Iuqm$,
so we may assume that $e > 1$.
Since $\Iuqm \subseteq I_x \cap I_y$
we only need to prove that $I_x^e \cap I_y^e \subseteq \Iuqm^e$.
Let $w$ be a monomial in $I_x^e \cap I_y^e$.
Then $w$ is a product of $e$ monomials in $I_x$
and it is also a product of $e$ monomials in $I_y$.
If one of the factors in the two products is identical,
then $w = w_0 w'$ for some $w_0 \in I_x \cap I_y = \Iuqm$
and some $w' \in I_x^{e-1} \cap I_y^{e-1}$,
and so by induction on $e$ we know that $w \in \Iuqm^e$.
Thus we may assume that the two factorizations of $w$
have no monomial-generator factor in common.
If in one of the factorizations the $e$ factors are in $I_1$
or more generally in $\Iuqm$,
then $w \in \Iuqm^e$,
and so we are done.
So in the writing of $w$ as an element of $I_1^e$
we may assume that at least one of the $e$ factors
is a monomial in $I_2 + I_3 x$,
and in the writing of $w$ as an element of $I_2^e$
we may assume that at least one of the $e$ factors
is a monomial in $I_3 y$.
In either factorization such a factor is of the form ${1 \over v} M_i Y_i$
for some $i \le k$ and $v = y$ in the former case
and $v = x$ in the latter case.
\proof
Let $x$ stand for $x_k$ or $y_{kl}$
and $y$ for $y_{jp} \not = x$.
By possibly renaming this also covers the case $y = y_{kl}$ when $x = x_k$.
We write $\Iuqm = I_1 + I_2 y + I_3 xy$,
where
$$
\eqalignno{
I_1 &= \Ioneqm Y_n + (M_i Y_i : i \in J_{\uq}, i > j), \cr
I_2 &= \cases{
{1 \over y} (M_i Y_i : i \in J_{\uq}, i \not = k, i \le j), & if $x = x_k$; \cr
{1 \over y} (M_i Y_i : i \in J_{\uq}, k < i \le j), & if $x = y_{kl}$,  \cr
} \cr
I_3 &= \cases{
{1 \over xy} (M_k Y_k), & if $x = x_k$; \cr
{1 \over xy} (M_i Y_i : i \in J_{\uq}, i \le k), & if $x = y_{kl}$. \cr
} \cr
}
$$
The generators of $I_1, I_2, I_3$ do not involve $x$ or $y$.
It suffices to prove that
$$
\Iuqm^e = (I_1 + I_2 + I_3 x)^e \cap (I_1 + I_2 y + I_3 y)^e,
$$
because $y$ is not in any prime ideal associated to the first intersectand
and $x$ is not in any prime ideal associated to the second intersectand.
Since for monomial ideals the intersection commutes with sums,
if $e = 1$,
the intersection equals
$I_1 + I_2 y + I_3 xy = \Iuqm$,
so we may assume that $e > 1$.
Since $\Iuqm^e$ is in the intersection,
it suffices to prove that any monomial~$w$ in the intersection
is also in $\Iuqm^e$.
We can write
$$
w = U U_{11} \cdots U_{1u_1} U_{21} \cdots U_{2u_2} U_{31} \cdots
U_{3u_3} x^{u_3}
= V V_{11} \cdots V_{1v_1} V_{21} \cdots V_{2v_2} V_{31} \cdots
V_{3v_3} y^{v_2 + v_3},
$$
where $U_{ir}, V_{is}$ are monomial generators of $I_i$,
$U, V$ are monomials,
and $u_1, u_2, u_3, v_1, v_2, v_3$ are non-negative integers
such that $u_1 + u_2 + u_3 = v_1 + v_2 + v_3 = e$.
We may assume that the greatest common divisor of $U$ and $V$ is~$1$.

By unique factorization,
$x^{u_3}$ divides $V$ and $y^{v_2+v_3}$ divides $U$.
If $u_3 \ge v_3$
then $w \in I_1^{v_1} (I_2y)^{v_2} (I_3 xy)^{v_3} \in \Iuqm^e$,
and if $v_2 + v_3 \ge u_2 + u_3$,
then
$w \in I_1^{u_1} (I_2y)^{u_2} (I_3 xy)^{u_3} \in \Iuqm^e$.
So we may assume that $u_3 < v_3$ and $v_2 + v_3 < u_2 + u_3$.

In particular,
$u_2 > v_2 + 1$,
and $u_2 + u_3$ and $v_3$ are positive.

If $a$ divides $U$ or $V$,
since $i \le n-1 < m$ and since $a \cdot a^m b^{m-i+1} = (a^{m+1}b) \in I_1$,
we may rewrite $w$ in the two forms as above
with strictly smaller $u_2+u_3$ or $v_3$
and correspondingly strictly larger $u_1$ or~$v_1$.
Thus we may assume that $a$ is not a factor of $U$ or of $V$.

If some $U_{1r}$ equals some $V_{1s}$
then by induction on $e$ we have that
$w/U_{1r} \in \Iuqm^{e-1}$
so that $w \in \Iuqm^e$.
So we may assume that
$\{V_{11}, \ldots, V_{1u_1}\} \cap \{U_{11}, \ldots, U_{1v_1}\} = \emptyset$.

Suppose that some $x_i$ with $i > j$ divides~$w$.
Then necessarily either $x_i$ divides $V$ and some $U_{1l}$
or else $x_i$ divides $U$ and some $V_{1l}$.
By \ref{lmreduction}~c),
$x_i (I_2 + I_3) \subseteq (M_i Y_i) \in~I_1$.
Thus we may rewrite $w$ so that $M_i Y_i$ appears as a $U_{1l}$
and as a $V_{1l'}$,
after which the two expressions for $w$ have a common factor in $I_1$,
and then as in the previous paragraph we are done by induction.
So we may assume that $x_i$ does not divide $w$ for all $i > j$.
Thus all $U_{1r}$ and $V_{1s}$ are in $\Ioneqm$.

If some $U_{1r}$ equals $a^{m+2}$ or $ab^{m+1}$,
then by \ref{lmreduction}
we may rewrite $w$ to increase $u_1$ by~$1$ and decrease $u_2 + u_3$ by~$1$.
Thus without loss of generality $U_{1r} \in (a^{m+1} b, b^{m+2})$.
Similarly,
by possibly decreasing~$v_3$
without loss of generality all $V_{1s} \in (a^{m+1} b, b^{m+2})$.
In particular,
either the product of the $U_{1r}$ equals $(a^{m+1}b)^{u_1}$
and the product of the $V_{1s}$ equals $b^{(m+2)v_1}$,
or else
the two products are respectively $b^{(m+2)u_1}$ and $(a^{m+1}b)^{v_1}$.
Then the $a$-degree of $w$ equals either
$(m+1) u_1 + m(u_2 + u_3) = m(v_2 + v_3)$
or
$m(u_2 + u_3) = (m+1) v_1 + m(v_2 + v_3)$.
Since $u_1, v_1 \le e \le k < n \le m$,
in the former case necessarily $u_1 = 0$
and in the latter case $v_1 = 0$.
In either case it follows that $u_2 + u_3 = v_2 + v_3$,
which contradicts the assumption
$u_2 + u_3 > v_2 + v_3$.
This proves the proposition.